\documentclass[]{article}

\usepackage[utf8]{inputenc}
\usepackage[T1]{fontenc}
\usepackage[hmarginratio=1:1,top=32mm,columnsep=20pt]{geometry} 
\usepackage[hang, small,labelfont=bf,up,textfont=it,up]{caption} 
\usepackage{booktabs} 
\usepackage{abstract} 

\usepackage{titlesec}
\titleformat{\section}[block]{\large\scshape\centering}{\thesection.}{1em}{} 
\titleformat{\subsection}[block]{\large}{\thesubsection.}{1em}{} 

\usepackage{amssymb, epsfig,amssymb, latexsym}
\usepackage{bm}
\usepackage{amsfonts,psfrag,amsmath,bbm,color,url}
\usepackage{graphicx}
\graphicspath{ {./img/} }
\usepackage{subfig}

\usepackage{amsfonts}
 \usepackage{amsmath}
 \usepackage{xcolor}
\usepackage{tikz}
\usepackage{pgfplots}
\usepackage{comment}
\usepackage{bm}
\usepackage{mathabx}
\DeclareMathAlphabet{\mathbbx}{U}{bbold}{m}{n}

\usepackage{cases}
\usepackage[output-decimal-marker={,}]{siunitx}
\usepackage{lineno,hyperref}
\definecolor{vargreen}{rgb}{0.0, 0.5, 0.0}
\definecolor{navyblue}{rgb}{0.0, 0.0, 0.5}
\newcommand{\green}[1]{{\color{black} #1}}
\newcommand{\blue}[1]{{\color{black} #1}}
\definecolor{mediumorchid}{rgb}{0.73, 0.33, 0.83}
\definecolor{crimson}{rgb}{0.86, 0.08, 0.24}
\definecolor{lightseagreen}{rgb}{0.13, 0.7, 0.67}
\definecolor{royalblue}{rgb}{0.25, 0.41, 0.88}
\definecolor{hotpink}{rgb}{1.0, 0.41, 0.71}
\definecolor{magenta}{rgb}{1.0, 0.0, 1.0}
\definecolor{goldenrod}{rgb}{0.85, 0.65, 0.13}

\definecolor{plum(traditional)}{rgb}{0.56, 0.27, 0.52}

\usepackage[export]{adjustbox}
\usepackage{mdframed}
\usepackage{array}
\usepackage{tikz}
\usetikzlibrary{shapes}
\usetikzlibrary{plotmarks}

\newtheorem{remark}{Remark}[section]

\def\PP{{{\rm l}\kern - .15em {\rm P} }}
\def\PN2{{\PP_{N}-\PP_{N-2}}}

\newcommand{\btau}{\boldsymbol{\tau}}

\newcommand{\ba}{\boldsymbol{a}}

\newcommand{\bb}{\boldsymbol{b}}

\newcommand{\bff}{\boldsymbol{f}}

\newcommand{\bu}{\boldsymbol{u}}

\newcommand{\red}[1]{{\color{black}#1}}

\newcommand{\deleted}[1]{{}}

\pgfplotsset{compat=1.17}
\usepackage[mode=buildnew]{standalone}
\title{\textbf{Hybrid Data-Driven Closure Strategies for Reduced Order Modeling}}

\date{ }
\pgfplotsset{compat=1.17}
\author{Anna Ivagnes  \\ \small SISSA, International School for Advanced Studies, \\ \small Mathematics Area, mathLab, Trieste, Italy. \\ \small  \href{mailto:aivagnes@sissa.it}{aivagnes@sissa.it} \normalsize \and Giovanni Stabile \\ \small SISSA, International School for Advanced Studies, \\ \small Mathematics Area, mathLab, Trieste, Italy. \\ \small  \href{mailto:gstabile@sissa.it}{gstabile@sissa.it} \normalsize  \and Andrea Mola  \\\small  Multi-scale Analysis of Materials Unit, \\ \small Scuola IMT Alti Studi, Lucca, Italy. \\ \small \href{mailto:andrea.mola@imtlucca.it}{andrea.mola@imtlucca.it} \normalsize \and Traian Iliescu  \\\small 
Department of Mathematics, Virginia Tech, \\ \small Blacksburg, VA, USA. \\ \small \href{mailto:iliescu@vt.edu}{iliescu@vt.edu} \normalsize \and Gianluigi Rozza  \\ \small  SISSA, International School for Advanced Studies, \\\small  Mathematics Area, mathLab, Trieste, Italy. \\ \small \href{mailto:grozza@sissa.it}{grozza@sissa.it}}

\begin{document}
\maketitle
\begin{abstract}
\noindent In this paper, we 
propose hybrid data-driven ROM closures for fluid flows.
These new ROM closures combine 
two fundamentally different strategies: 
(i) purely data-driven ROM closures, both for the velocity and the pressure; and
(ii) physically based, eddy viscosity  data-driven closures, which model the energy transfer in the system. 
The first 
strategy 
consists in the addition of closure/correction terms 
to the governing equations, which are built from the available data. 
The second 
strategy includes turbulence modeling by adding eddy viscosity terms, which are determined by using machine learning techniques.
The 
two strategies 
are combined for the first time in this paper to investigate 
a two-dimensional flow past a circular cylinder at $Re=\num{50000}$. 
Our numerical results show that the hybrid data-driven ROM is more accurate than both the purely data-driven ROM and the eddy viscosity ROM.
\end{abstract}



\section{Introduction}



Full order models (FOMs) are computational models that are obtained by using classical numerical methods, such as the finite element or the finite volume methods.
The FOM computational cost is prohibitively high for many important engineering, environmental, and biomedical applications that require repeated turbulent flow simulations, such as \cite{stabile2020efficient}. 

Reduced Order Models (ROMs), presented in \cite{degruyter1, degruyter2, degruyter3,rozza2008reduced, rozza2013reduced,aromabook} and applied to finite-volume schemes in \cite{haasdonk2008reduced, carlberg2018conservative}, can reduce the FOM computational cost by orders of magnitude.
ROMs can be built by using different methodologies, such as least-squares Petrov--Galerkin projection \cite{carlberg2011efficient} or Galerkin projection \cite{noack1994low, bergmann2009enablers, kunisch2002galerkin, Azaiez2017}.
In this paper, we consider the class of Galerkin-ROMs (G-ROMs), which are ROMs constructed by using a Galerkin method. The standard methodology consists in the construction of a set of basis functions (modes) for the velocity and pressure fields, $\{ \varphi_{1}, \ldots, \varphi_{r} \}$ and $\{ \chi_{1}, \ldots, \chi_{q} \}$, respectively, such that the unknown solution can be approximated as a linear combination of these basis functions: 
$u(x,t) = \sum_{i=1}^{r} a_i(t) \, \varphi_i(x)$, $p(x,t) = \sum_{i=1}^{q} b_i(t) \, \chi_i(x)$. The 
G-ROM consists in the resolution of a dynamical system, which is obtained by projecting the governing equations onto the space spanned by these modes.
The ROM basis is a data-driven basis, which is built by 
available FOM snapshots.
In this paper, we use the proper orthogonal decomposition (POD)~\cite{HLB96} to build the ROM basis.
We note, however, that other approaches can be used (see, e.g.,~\cite{hesthaven2015certified,quarteroni2015reduced}).
The resulting Galerkin ROM (G-ROM) is a system of equations in which the unknowns are the coefficients $(a_i)_{i=1}^{N_u}$ and $(b_i)_{i=1}^{N_p}$ appearing in the reduced fields expressions.
For example, for the incompressible Navier-Stokes equations, the G-ROM can be written as follows:
\begin{equation}
    \begin{cases}
        \ba_t 
        = \bff(\ba,\bb), & \\
        {\bf h}(\ba,\bb)
        = {\bf 0} , & 
    \end{cases} 
    \label{g-rom}
\end{equation}
where $\ba(t) := (a_1(t), \ldots, a_r(t))^T$ and $\bb(t) := (b_1(t), \ldots, b_q(t))^T$, and $\bff$ and $\bf h$ denote the ROM operators. 

ROMs are often built 
by using a few 
modes. 
In the {\it marginally-resolved} modal regime, where the number of ROM basis functions 
allows a moderately accurate representation of 
the main features of the underlying dynamics,
G-ROMs can yield acceptable approximations. 
However, in the case of turbulent flow applications, hundreds or even thousands of ROM modes are often required to obtain a good approximation of the flow dynamics. 
Therefore, 
in those cases (i.e., in the {\it under-resolved} regime),
low-dimensional G-ROMs generally yield inaccurate approximations.

Different stabilization approaches have been studied in the past years in order to avoid the stability issues and to increase the ROM accuracy \cite{grimberg2020stability, bergmann2009enablers, iollo2000stability, kalashnikova2014construction}.

In this paper, in the marginally-resolved and, especially, under-resolved regimes, in
order to increase the accuracy of the standard G-ROM \eqref{g-rom} while maintaining an acceptable 
computational cost, the ROM approach is integrated with data-driven techniques. One popular approach is to add 
\emph{closure/correction} terms 
to the reduced formulation. 
To this end, one can use a 
\emph{purely data-driven} 
approach, which exploits only the 
available FOM data to build new terms 
that include the contribution of the modes neglected in the ROM formulation.
The new reduced system with 
purely data-driven 
terms is the following:
\begin{equation}
    \begin{cases}
        \ba_t 
        = \bff(\ba,\bb)
        + \btau_{\bu}(\ba, \bb), & \\
        {\bf h}(\ba,\bb)
        + \btau_{p}(\ba, \bb)
        = {\bf 0} . & 
    \end{cases}
    \label{eqn:p-g-rom+closure}  
\end{equation}
The terms $\btau_{u}$ and $\btau_p$ are evaluated by solving a least square problem 
that minimizes the difference between the model form of the correction terms and the exact correction terms, which are evaluated by using the available FOM data 
\cite{paper1,xie2018data, mohebujjaman2019physically, mou2021data}. The technique used to compute the correction terms takes inspiration from the operator inference approach \cite{peherstorfer2016data, wilcoxNLROM2019, benner2020operator, karasozen2022intrusive}.
In \cite{paper1}, 
the purely data-driven approach was applied for the first time to the pressure Poisson formulation (PPE-ROM), leading to the introduction of 
the novel pressure correction term, $\btau_p$. This 
strategy has significantly increased the accuracy 
of the pressure field approximation.


To construct the ROM correction terms, a different  strategy is the \emph{physically-based data-driven} approach.
This strategy consists of two steps:
In the first step, one postulates a physical model form for the ROM correction term.
For example, to approximate turbulent flows, one can postulate an eddy viscosity model.
In the second step of the physically-based data-driven strategy, one solves a least squares problem to determine the optimal parameters in the ROM correction term, e.g., the eddy viscosity coefficient.
The physically-based data-driven strategy was used in~\cite{hijazi2020data}, where the FOM utilizes a Reynolds averaged Navier-Stokes (RANS) approach \cite{reynolds1895iv} and an eddy viscosity model is leveraged to close the FOM system \cite{jones1972prediction, kolmogorov1941equations, spalding1974numerical}. 
To enforce FOM-ROM model consistency, a turbulence model was also introduced at a ROM level in~\cite{rebollo2017certified} and in \cite{hijazi2020data}, where an eddy viscosity reduced order model (EV-ROM) is built.
In the EV-ROM, 
the eddy viscosity reduced field is expressed as a linear combination of eddy viscosity modes $\nu_t(x,t)=\sum_{i=1}^{N_{\nu_t}} g_i(t) \eta_i(x)$, where $\{\eta_1, \ldots, \eta_{N_{\nu_t}}\}$ \blue{are the eddy viscosity modes}. The model introduced in \cite{hijazi2020data} can be expressed as follows:
\begin{equation}
    \begin{cases}
        \ba_t 
        = \bff(\ba,\bb,\mathbf{g}), & \\
        {\bf h}(\ba,\bb,\mathbf{g})
        = {\bf 0} , & 
    \end{cases} 
    \label{turb-g-rom}
\end{equation}
The vector of eddy viscosity coefficients $\mathbf{g}:=(g_1(t), \ldots, g_{N_{\nu_t}})$ can be obtained from the velocity vector of coefficients $\ba$ 
by using either interpolation or regression techniques.
\medskip

The main aim of this paper is to propose 
a novel \emph{hybrid data-driven} ROM that is constructed by 
combining the two fundamentally different strategies outlined above: purely data-driven and physically-based data-driven techniques. 
The new hybrid data-driven ROM 
can be written as follows:
\begin{equation}
    \begin{cases}
        \ba_t 
        = \bff(\ba,\bb,\mathbf{g})
        + \btau_{\bu}(\ba, \bb), & \\
        {\bf h}(\ba,\bb)
        + \btau_{p}(\ba, \bb,\mathbf{g})
        = {\bf 0} . & 
    \end{cases}
    \label{turb-p-g-rom+closure}  
\end{equation}
\emph{Machine learning} techniques are exploited to find the eddy viscosity coefficients in~\eqref{turb-p-g-rom+closure}, training a fully-connected neural network from the full order data.

To our knowledge, this is the first time that a purely data-driven strategy for both the velocity and the pressure, and a physically-based eddy viscosity data-driven strategy are combined.
A similar 
approach has been used in~\cite{mou2021data-phd}.
We note, however, that the approach in~\cite{mou2021data-phd} is different from the approach used in the present study in two significant aspects:
first, a data-driven correction term is used in this paper, but not in~\cite{mou2021data-phd};
second, the physically-based data-driven correction terms are different: regression based mixing length in~\cite{mou2021data-phd}, and machine learning based eddy viscosity in the current study.

In this paper, the 
new hybrid 
data-driven ROM ~\eqref{turb-p-g-rom+closure} is compared to (i) the purely data-driven ROM \eqref{eqn:p-g-rom+closure} proposed in \cite{paper1}, and 
(ii) the physically-based data-driven ROM \eqref{turb-g-rom} proposed in \cite{hijazi2020data}. 
The three 
ROMs are tested and compared on the classical case study of a turbulent flow around a cylinder. The numerical results 
show that the 
new data-driven ROM~\eqref{turb-p-g-rom+closure} yields more accurate velocity and pressure approximations than both the purely data-driven ROM \eqref{eqn:p-g-rom+closure} and the physically-based data-driven ROM \eqref{turb-g-rom}.

These results support the following two conclusions, which are similar to those in~\cite{mou2021data-phd}:
adding a physically-based (eddy viscosity) correction term improves the accuracy of the purely data-driven ROM.
Furthermore, adding purely data-driven velocity and pressure correction terms improves the accuracy of the physically-based data-driven ROM.
Thus, our numerical investigation suggests that a hybrid data-driven ROM closure strategy is more accurate  than both a purely data-driven closure and a physically-based data-driven closure.

The rest of the paper is organized as follows: Section
\ref{sec:fom} is dedicated to a brief overview of the 
FOM used and implemented in the C++ open source software \emph{OpenFOAM}. 
In section \ref{sec:g-rom}, the standard G-ROM framework is recalled, specifying the two approaches used in 
the numerical investigation 
(the supremizer and 
pressure Poisson approaches). Sections \ref{dataROM} and \ref{evROM} summarize the purely data-driven ROMs and physically-based data-driven ROMs, respectively.
Section \ref{combined_dataROM} is dedicated to the presentation of the new hybrid data-driven ROM, 
whereas 
section \ref{results} 
presents the results of the numerical simulations in terms of relative errors of the pressure and velocity reduced fields with respect to the full order results.
Finally, section~\ref{sec:conclusions} presents the conclusions of our study and outlines future research directions.

\section{Full Order Model (FOM)} 
\label{sec:fom}
As a mathematical model, we use 
the Navier-Stokes Equations (NSE) for incompressible flows. 
The fluid domain is $\Omega \in \mathbb{R}^d$ with $d=
2$ or $3$, $\Gamma$ 
the domain's boundary, 
$t \in [0,T]$ 
the time, $\mathbf{u}=\mathbf{u}(\mathbf{x},t)$ 
the flow velocity vector field, 
$p=p(\mathbf{x},t)$ 
the pressure scalar field normalized by the fluid density, and $\nu$ 
the fluid kinematic viscosity. 
The strong form of the incompressible NSE is the following:
\begin{subnumcases}{\label{NSE}}
\frac{\partial \mathbf{u}}{\partial t}=-\nabla \cdot (\mathbf{u} \otimes \mathbf{u})+\nabla \cdot \nu \left(\nabla \mathbf{u} + (\nabla \mathbf{u})^T \right)-\nabla p & in $\Omega \times [0,T]\, , $ \label{mom_NSE}\\
\nabla \cdot \mathbf{u} = \mathbf{0} & in $\Omega \times [0,T]\, ,$ \label{cont_NSE}\\
+ \text{ boundary conditions }& on  $\Gamma \times [0,T]\, ,$ \label{bound_NSE}\\
+ \text{ initial conditions }& in  $(\Omega,0)\, .$ \label{init_NSE}
\end{subnumcases}

The full order solutions of \ref{NSE} are computed by 
using the software \emph{OpenFOAM}, which employs a finite volume discretization of the system \cite{moukalled2016finite, jasak1996error}.


The full order problem 
uses the U-RANS (
unsteady RANS equations) approach. 
This approach is based on the 
Reynolds decomposition, 
proposed in \cite{reynolds1895iv}, which decomposes each flow field in its mean and its fluctuating part. The RANS equations are obtained by taking the time-average of the NSE in \eqref{NSE}. In the resulting system, 
it is important to include a closure model for the 
Reynolds stress tensor. 
The closure model 
considered in this paper is an 
eddy viscosity 
model 
that is based on the Boussinesq hypothesis. Specifically, 
the SST $\kappa-\omega$ model is used to close the RANS system. This model is based on the inclusion of two transport equations 
that describe two additional variables: the kinetic energy, $\kappa$, and the specific turbulent dissipation rate, $\omega$. The 
SST $\kappa-\omega$ model is presented in \cite{kolmogorov1941equations} in 
the standard formulation, and in \cite{menter1994two} in the SST formulation. The extended RANS model including the SST $\kappa-\omega$ equations can be found in \cite{hijazi2020data}.

\section{The POD-Galerkin ROMs}
    \label{sec:g-rom}
In this section, a brief overview of the standard POD-Galerkin ROM techniques is provided.
\medskip

After the offline stage is performed, all the FOM snapshots, i.e., the FOM solutions for different time instants $\{t_j\}_{j=1}^{N_T}$, 
are collected. 
The 
POD is then applied on the 
full order snapshot matrices:
\begin{equation*}
\mathcal{S}_u=\{\mathbf{u}(\mathbf{x},t_1),...,\mathbf{u}(\mathbf{x},t_{N_T})\}  \in \mathbb{R}^{N_u^h \times N_T}, \quad \mathcal{S}_p=\{p(\mathbf{x},t_1),...,p(\mathbf{x},t_{N_T})\}  \in \mathbb{R}^{N_p^h \times N_T},
\end{equation*}
where $N_u^h$ and $N_p^h$ are the numbers 
of degrees of freedom for the velocity and pressure fields, respectively.

The ROM construction uses the following
different stabilization approaches for the velocity-pressure coupling 
\cite{paper1}:
\begin{itemize}
    \item[(1)] the SUP-ROM approach, in which additional \emph{supremizer} modes for the velocity space are introduced in order to fulfill the \emph{inf-sup} condition \cite{rozza2007stability, ballarin2015supremizer, stabile2018finite, ali2020stabilized};
    \item[(2)] the PPE-ROM approach, in which the pressure Poisson equation replaces the continuity equation \cite{akhtar2009stability,Stabile2017CAIM,stabile2018finite, noack2005need}.
\end{itemize}

In the two formulations, the velocity and pressure POD spaces are assembled as follows:
\begin{equation}
\mathbb{V}^u_{\text{POD}}=\mbox{span}\{[\boldsymbol{\phi}_i]_{i=1}^{N_u}, [\mathbf{s}(\chi_i)]_{i=1}^{N_{sup}}
\},\quad 
\mathbb{V}^p_{\text{POD}}=\mbox{span}\{[\chi_i]_{i=1}^{N_p}\},
\end{equation}
where $N_u \ll N_u^h$ and $N_p\ll N_p^h$, and $[\boldsymbol{\phi}_i]_{i=1}^{N_u}$ and $[\chi_i]_{i=1}^{N_p}$ are the velocity and pressure POD modes, respectively.
\blue{
The supremizer modes $(\mathbf{s}_i)_{i=1}^{N_{sup}}=\mathbf{s}(\chi_i)_{i=1}^{N_{sup}}$ are additional modes introduced in the SUP-ROM approach in order to fulfill the \emph{inf-sup} condition \cite{stabile2018finite}.

For each pressure basis function, the corresponding supremizer element can be found by solving the following problem:
\begin{equation}
\begin{cases}
 \Delta \mathbf{s}_i = - \nabla p_i \text{ in } \Omega , \\
 \mathbf{s}_i = 0 \text{ on } \partial \Omega .
\end{cases}
\label{sup_condition}
\end{equation}

The velocity POD space can be enriched either with an \emph{exact} or with an \emph{approximated} approach \cite{ballarin2015supremizer}. In the exact approach, the problem \eqref{sup_condition} is solved for each pressure basis function $\chi_i$ and each solution is added to the velocity space. In the approximated approach, the problem \eqref{sup_condition} is solved for each pressure
snapshot $p(\mathbf{x}, t_i)_{i=1}^{N_T}$, which yields the following supremizer snapshot matrix:
\begin{equation*}
\mathcal{S}_{sup}=\{\mathbf{s}(\mathbf{x},t_1),...,\mathbf{s}(\mathbf{x},t_{N_T})\}  \in \mathbb{R}^{N_u^h \times N_T}.
\end{equation*}
A POD modal decomposition is then applied to the snapshot matrix in order to obtain the supremizer POD modes $(\boldsymbol{\eta}_i)_{i=1}^{N_{sup}}$ \cite{stabile2018finite}.
In this paper, we adopt the approximated procedure, since it significantly reduces the computational cost of the offline phase.
}
The reduced velocity and pressure fields are expressed as follows:
\begin{equation}
\mathbf{u}(\mathbf{x},t) \approx \mathbf{u}_r(\mathbf{x},t)=\sum_{i=1}^{N_u+N_{sup}} a_i(t)\boldsymbol{\phi}_i(\mathbf{x}), \quad
p(\mathbf{x},t)\approx p_r(\mathbf{x},t)=\sum_{i=1}^{N_p} b_i(t)\chi_i(\mathbf{x}),
\label{appfields}
\end{equation}
where $N_{sup}=0$ for the PPE-ROM approach. 


In the remainder of this section, we outline the SUP-ROM (section~\ref{sec:sup-rom}) and PPE-ROM (section~\ref{sec:ppe-rom}) approaches.

\subsection{SUP-ROM}
    \label{sec:sup-rom}

In the supremizer 
approach, performing a Galerkin projection of the momentum equation \eqref{mom_NSE} onto the velocity modes, and of the continuity equation \eqref{cont_NSE} onto the pressure modes, the following dynamical system is obtained:
\begin{equation}
    \begin{cases}
    \mathbf{M} \dot{\mathbf{a}}=\nu(\mathbf{B}+\mathbf{B_T})\mathbf{a}-\mathbf{a}^T \mathbf{C} \mathbf{a}-\mathbf{H}\mathbf{b}+\tau \left( \sum_{k=1}^{N_{\text{BC}}}(U_{\text{BC},k}\mathbf{D}^k-\mathbf{E}^k \mathbf{a})\right)\, ,\\
    \mathbf{P}\mathbf{a}=\mathbf{0}\, ,
    \end{cases}
    \label{reduced_system}
\end{equation}
where $\mathbf{a}$ and $\mathbf{b}$ are the vectors of the coefficients associated to the velocity and pressure modes, respectively.

The matrices appearing in the system are defined 
as follows:
\[
\begin{split}
&(\mathbf{M})_{ij}=(\boldsymbol{\phi}_i,\boldsymbol{\phi}_j)_{L^2(\Omega)}, \quad (\mathbf{P})_{ij}=(\chi_i,\nabla \cdot \boldsymbol{\phi}_j)_{L^2(\Omega)}\, ,\quad (\mathbf{B})_{ij}=(\boldsymbol{\phi}_i,\nabla \cdot \nabla \boldsymbol{\phi}_j)_{L^2(\Omega)}, \\
&(\mathbf{B_T})_{ij}=(\boldsymbol{\phi}_i,\nabla \cdot (\nabla \boldsymbol{\phi}_j)^T)_{L^2(\Omega)},\quad (\mathbf{C})_{ijk}=(\boldsymbol{\phi}_i,\nabla \cdot (\boldsymbol{\phi}_j \otimes \boldsymbol{\phi}_k))_{L^2(\Omega)}, \quad (\mathbf{H})_{ij}=(\boldsymbol{\phi}_i,\nabla \chi_j)_{L^2(\Omega)}\, .
\end{split}
\]
\green{\remark{We highlight that the POD modes, as they are linear combination of snapshots, are piece-wise constant functions. However, it is possible to compute second order operators using the same FV procedure used at the FOM level. We remind to \cite{moukalled2016finite} for details.}}

The term $\tau \left( \sum_{k=1}^{N_{\text{BC}}}(U_{\text{BC},k}\mathbf{D}^k-\mathbf{E}^k \mathbf{a})\right)$ in \eqref{reduced_system} is a term used to enforce the Dirichlet boundary conditions in the reduced order model \cite{hijazi2020data, star2019extension}. 
\blue{
We call $\Gamma_D$ the Dirichlet boundary, which is composed by different parts, $\Gamma_{D_1}$, $\Gamma_{D_2}$, .., $\Gamma_{D_K}$.
}

In particular, $N_{\text{BC}}$ is the number of velocity boundary conditions \blue{we would like to impose on the parts} of the Dirichlet boundary \blue{where velocity has at least one non-zero component}; $U_{\text{BC},k}$ is the velocity non-zero component at the $k$-th part of the Dirichlet boundary.

$\tau$ is a penalization factor \blue{which is tuned by a sensitivity analysis on the specific problem comsidered \cite{hijazi2020data, star2019extension}. In general, bigger values of the penalization factor lead to a stronger enforcement of the boundary conditions.}

The matrices $\mathbf{E}^k$ and vectors $\mathbf{D}^k$ are defined as follows:
\[(\mathbf{E}^k)_{ij}=(\boldsymbol{\phi}_i, \boldsymbol{\phi}_j)_{L^2(\Gamma_{D_k})}, \quad (\mathbf{D}^k)_{i}=(\boldsymbol{\phi}_i)_{L^2(\Gamma_{D_k})}, \text{ for all }k=1,...,N_{\text{BC}}.\]

\subsection{PPE-ROM}
    \label{sec:ppe-rom}

In the pressure Poisson approach, the continuity equation \eqref{reduced_system} is replaced by the Poisson equation for pressure, obtained by taking the divergence of the momentum equation. 
The dynamical system is obtained by projecting the momentum equation and the pressure Poisson equation on the velocity and pressure modes, respectively:
\begin{equation}
    \begin{cases}
    \mathbf{M} \dot{\mathbf{a}}=\nu(\mathbf{B}+\mathbf{B_T})\mathbf{a}-\mathbf{a}^T \mathbf{C} \mathbf{a}-\mathbf{H}\mathbf{b}+\tau \left( \sum_{k=1}^{N_{\text{BC}}}(U_{BC,k}\mathbf{D}^k-\mathbf{E}^k \mathbf{a})\right)\, ,\\
    \mathbf{D}\mathbf{b}+ \mathbf{a}^T \mathbf{G} \mathbf{a} - \nu \mathbf{N} \mathbf{a}- \mathbf{L}=\mathbf{0}\, .
    \end{cases}
    \label{reduced_systemPPE}
\end{equation}

The matrices $\mathbf{M}, \, \mathbf{B}, \, \mathbf{B_T}, \, \mathbf{C}, \, \mathbf{H}, \, \mathbf{D}^k$, and $\mathbf{E}^k$ in system \eqref{reduced_systemPPE} are the same as the corresponding matrices in~\eqref{reduced_system}.
The additional matrices 
are defined as follows:
\[
\begin{split}
&(\mathbf{D})_{ij}=(\nabla \chi_i,\nabla \chi_j)_{L^2(\Omega)}, \quad 
(\mathbf{G})_{ijk}=(\nabla \chi_i,\nabla \cdot (\boldsymbol{\phi}_j \otimes \boldsymbol{\phi}_k))_{L^2(\Omega)}, \\ &(\mathbf{N})_{ij}=(\mathbf{n} \times \nabla \chi_i,\nabla \boldsymbol{\phi}_j)_\blue{{L^2(\Gamma)}}, \quad (\mathbf{L})_{ij}=(\chi_i,\mathbf{n} \cdot \boldsymbol{R}_t)_\blue{{L^2(\Gamma)}}\, \blue{,}
\end{split}
\]
\blue{where vector $\mathbf{n}$ is the normal unitary vector to the domain boundary $\Gamma$.}

\section{Data-driven VMS-ROMs}
\label{dataROM}
We consider the \emph{data-driven variational multiscale ROM} (\emph{DD-VMS-ROM}) framework~\cite{paper1}, in which different 
correction/closure terms are introduced in the formulations of SUP-ROM
\eqref{reduced_system} and PPE-ROM \eqref{reduced_systemPPE}. 
The DD-VMS-SUP-ROM approach (section~\ref{datasup}) includes only velocity correction terms, whereas the 
DD-VMS-PPE-ROM approach (section~\ref{datappe}) includes both pressure and velocity correction terms.

\blue{
For the sake of simplicity, in this section we will consider the following notation for the coefficient vectors for velocity and pressure:
\[
\mathbf{a}=(a_i)_{i=1}^r, \quad \mathbf{b}=(b_i)_{i=1}^q,
\]
where $r$ is the reduced number of modes for velocity ($r=N_u$ in the PPE-ROM approach, $r=N_u+N_{sup}$ in the SUP-ROM approach), $q$ is the reduced number of modes for pressure ($q=N_p$). We will also consider $r_{tot}=r+q$.
}

\subsection{DD-VMS-SUP-ROM}
\label{datasup}
The reduced formulation for the supremizer approach adopted in 
\cite{paper1} is the following: 
\begin{equation}
\begin{cases}
\mathbf{M} \dot{\mathbf{a}}=\nu(\mathbf{B}+\mathbf{B_T})\mathbf{a}-\mathbf{a}^T \mathbf{C} \mathbf{a}-\mathbf{H}\mathbf{b}+\tau \left( \sum_{k=1}^{N_{BC}}(U_{BC,k}\mathbf{D}^k-\mathbf{E}^k \mathbf{a})\right)+\boldsymbol{\tau}^u(\ba)\, ,\\
    \mathbf{P}\mathbf{a} =\mathbf{0}\, .
\end{cases}
\label{new_reduced_system}
\end{equation}
The velocity correction term appearing in \eqref{new_reduced_system} is modeled as follows:
\begin{equation}
\boldsymbol{\tau}^u (\mathbf{a})=\tilde{A} \mathbf{a}+\mathbf{a}^T \tilde{B} \mathbf{a},
\label{ansatz1}
\end{equation}
where $\tilde{A}$ is a matrix and $\tilde{B}$ is a three-dimensional tensor.
The 
operators $\tilde{A}$ and $\tilde{B}$ 
are constructed by solving an optimization problem that minimizes the difference 
between the exact correction term,  $\boldsymbol{\tau}_u^{\text{exact}}$, and 
a correction term ansatz, $\boldsymbol{\tau}_u^{\text{ansatz}}$:
\begin{equation}
\min_{\substack{\tilde{A} \in \mathbb{R}^{r \times r},\\ \tilde{B} \in \mathbb{R}^{r \times r \times r}}}{\sum_{j=1}^M || \boldsymbol{\tau}_u^{\text{exact}}(t_j)-\boldsymbol{\tau}_u^{\text{ansatz}}(t_j)||_{L^2(\Omega)}^2}\, ,
\label{opt_problem}
\end{equation}
where $M$ time instances are considered to build the correction term and the term $\boldsymbol{\tau}^{\text{exact}}(t_j)$ is computed starting from the snapshots coefficient vectors $\mathbf{a}_d^{snap}(t_j)$ and $\mathbf{a}_r^{snap}(t_j)$, defined as follows:
\[
\begin{split}
&a_{d_i}^{snap}(t_j)=\left(\mathbf{u}_d(t_j), \boldsymbol{\phi}_i\right)_{L^2(\Omega)} \quad \forall i=1,...,d\gg r.\\
&\blue{
a_{r_i}^{snap}(t_j)=\left(\mathbf{u}_r(t_j), \boldsymbol{\phi}_i\right)_{L^2(\Omega)} \quad \forall i=1,...,r.}
\end{split}\]
\red{
The value $d$ is the number of modes used to build the exact correction term. \remark{We remark that the value of $d$ can be chosen as the rank of the snapshots matrix. However, in order to decrease the computational cost of the offline stage, we chose a value which is smaller than the rank of the snapshots matrix, but considerably bigger than the reduced number of modes considered, as pointed out in \cite{xie2018data}.}}

The exact correction term for velocity is evaluated as follows:
\[ \boldsymbol{\tau}\blue{_u}^{\text{exact}}(t_j)=\left(- \overline{(\mathbf{a}_d^{snap}(t_j))^T\mathbf{C_d}\mathbf{a}_d^{snap}(t_j)}^r\right) -\left(-(\mathbf{a}_r^{snap}(t_j))^T\mathbf{C}\mathbf{a}_r^{snap}(t_j) \right)\, , \]
where the tensor $\mathbf{C_d} \in \mathbb{R}^{d \times d \times d}$ is defined in the following way:
\[
\mathbf{C_d}_{\,ijk}=\left(\boldsymbol{\phi}_i, \nabla \cdot (\boldsymbol{\phi}_j \otimes \boldsymbol{\phi}_k) \right)_{L^2(\Omega)}.
\]

The 
correction term ansatz is evaluated as in \eqref{ansatz1}, but starting from $\mathbf{a}_r^{snap}(t_j)$ at each time step $j$:
\begin{equation}
\boldsymbol{\tau}\blue{_u}^{\text{ansatz}}(t_j)=\tilde{A}\mathbf{a}_r^{snap}(t_j)+(\mathbf{a}_r^{snap}(t_j))^T\tilde{B}\mathbf{a}_r^{snap}(t_j).
\end{equation}
The optimization problem \eqref{opt_problem} is rewritten as a least squares problem following a procedure similar to that used in \cite{peherstorfer2016data} and presented in detail in \cite{paper1} and \cite{ivagnes2021data}. 

A different way to find $\tilde{A}$ and $\tilde{B}$, introduced in \cite{mohebujjaman2019physically}, is that of solving a constrained least squares problem, inheriting the physical properties of the exact tensors:
\begin{equation}
    \min_{\substack{\tilde{A} \in \mathbb{R}^{r \times r},\\ \tilde{B} \in \mathbb{R}^{r \times r \times r},\\ \mathbf{a}^T  \tilde{A} \mathbf{a} \leq 0,\\\mathbf{a}^T(\mathbf{a}^T \tilde{B} \mathbf{a})=0}}{\sum_{j=1}^M || \boldsymbol{\tau}\blue{_u}^{\text{exact}}(t_j)-\boldsymbol{\tau}\blue{_u}^{\text{ansatz}}(t_j)||_{L^2(\Omega)}^2}\, .
    \label{opt_problem2}
\end{equation}
As shown in 
\cite{paper1,mohebujjaman2019physically}, the constrained method 
can yield better results than the unconstrained 
method 
in the 
marginally-resolved
regime.

\subsection{DD-VMS-PPE-ROM}
\label{datappe}
In this section, the pressure data-driven model developed in \cite{paper1} is briefly recalled. The continuity equation at the reduced level is replaced by the pressure Poisson equation. This formulation allows for the introduction of novel pressure correction terms in the reduced system: 
\begin{equation}
\begin{cases}
\mathbf{M} \dot{\mathbf{a}}=\nu(\mathbf{B}+\mathbf{B_T})\mathbf{a}-\mathbf{a}^T \mathbf{C} \mathbf{a}-\mathbf{H}\mathbf{b}+\tau \left( \sum_{k=1}^{N_{BC}}(U_{BC,k}\mathbf{D}^k-\mathbf{E}^k \mathbf{a})\right)+\boldsymbol{\tau}^u(\ba, \bb)\,,\\
    \mathbf{D}\mathbf{b}+\mathbf{a}^T \mathbf{G} \mathbf{a} - \nu \mathbf{N} \mathbf{a} - \mathbf{L}+\boldsymbol{\tau}^{p}(\ba, \bb)=\mathbf{0}\,.
\end{cases}
\label{new_reduced_systemPPE}
\end{equation}
In 
\cite{paper1,ivagnes2021data},
different ansatzes are analyzed for the correction terms $\boldsymbol{\tau}^u$ and $\boldsymbol{\tau}^p$. The formulation which is chosen in this paper is the one in which a unique least squares problem is solved to find all the purely data-driven terms.

As the exact correction term in the optimization problem, we consider the following:
\begin{equation}
    \boldsymbol{\tau}_{\text{tot}}^{\text{exact}}(t_j) = \left(\boldsymbol{\tau}_u^{\text{exact}}(t_j),\boldsymbol{\tau}_p^{\text{exact}}(t_j)\right) \, \forall j=1,...,M\,.
\end{equation}
The exact correction corresponding to $\boldsymbol{\tau}_p$ includes the contribution of two different terms:
\[
\begin{split}
 &\boldsymbol{\tau}_p^{\text{exact}}(t_j)= \boldsymbol{\tau}_D^{\text{exact}}(t_j)+\boldsymbol{\tau}_G^{\text{exact}}(t_j) =\\&= \left( \overline{\mathbf{D_d}\mathbf{b}_d^{\text{snap}}(t_j)}^q \right)-\mathbf{D} \mathbf{b}_q^{\text{snap}}(t_j)+ \overline{(\ba_d^{\text{snap}}(t_j))^T \mathbf{G_d} \ba_d^{\text{snap}}(t_j)}^r -(\ba_r^{\text{snap}}(t_j))^T\mathbf{G} \ba_r^{\text{snap}}(t_j)\, ,   
\end{split}\]
where the matrix $\mathbf{D_d}$ and tensor $\mathbf{G_d}$ are defined as follows:
\[
\mathbf{D_d}_{ij}=(\nabla \chi_i,\nabla \chi_j)_{L^2(\Omega)}, \quad \mathbf{G_d}_{ijk} = (\nabla \chi_i, \nabla \cdot (\boldsymbol{\phi}_j \otimes \boldsymbol{\phi}_k))_{L^2(\Omega)}, \, i,j,k = 1,\ldots ,d\, .
\]
The ansatz considered in this paper is the following:
\begin{equation}
    \boldsymbol{\tau}_{\text{tot}}^{\text{ansatz}}(t_j)=\tilde{J}_A \mathbf{ab}^{\text{snap}} (t_j)+ (\mathbf{ab}^{\text{snap}} (t_j))^T \tilde{J}_B \mathbf{ab}^{\text{snap}} (t_j)\,,
\end{equation}
where the matrices $\tilde{J}_A \in \mathbb{R}^{(r+q) \times (r+q)}$ and $\tilde{J}_B \in \mathbb{R}^{(r+q) \times (r+q) \times (r+q)}$, and the vector $\mathbf{ab}^{\text{snap}} (t_j)=(\ba_r^{\text{snap}} (t_j), \bb_q^{\text{snap}} (t_j))$ $\in \mathbb{R}^{r+q}$. The final correction which is inserted in the reduced system \eqref{new_reduced_systemPPE} is divided 
into two vectors:
\begin{equation}
\tilde{J}_A \mathbf{ab} + \mathbf{ab}^T \tilde{J}_B \mathbf{ab}=\left( \boldsymbol{\tau}\blue{_u}, \boldsymbol{\tau}\blue{_p}\right), \text{ where }\boldsymbol{\tau}\blue{_u} \in \mathbb{R}^{N_u}, \boldsymbol{\tau}\blue{_p} \in \mathbb{R}^{N_p}\,.
\end{equation}
Finally, the following optimization problem is solved:
\begin{equation}
    \min_{\substack{\tilde{J}_A \in \mathbb{R}^{r_{\text{tot}} \times r_{\text{tot}}}; \\ \tilde{J}_B \in \mathbb{R}^{r_{\text{tot}} \times r_{\text{tot}} \times r_{\text{tot}}}}}{\sum_{j=1}^M || \boldsymbol{\tau}_{\text{tot}}^{\text{exact}}(t_j)-\boldsymbol{\tau}_{\text{tot}}^{\text{ansatz}}(t_j)||_{L^2(\Omega)}^2}\,.
\label{opt_problem_abC}
\end{equation}

\section{Eddy Viscosity ROMs}
\label{evROM}

The 
ROMs described in section \ref{dataROM} do not include 
turbulence treatment inside the formulation.
Different models have been used in fluid applications to simulate the turbulent behaviour. In this paper, 
turbulence is included at the full order level by using a particular type of \emph{eddy viscosity} model, the SST $\kappa-\omega$ model, which 
adds to the RANS equations the transport equations for $\kappa$ and $\omega$. 

At the reduced order level, an approximation of the eddy viscosity terms can be included in equations by introducing a reduced order version of the eddy viscosity \cite{hijazi2020data}, as follows:
\begin{equation*}
    \nu_t(\mathbf{x},t) \approx \nu_{t_r} (\mathbf{x},t)=\sum_{i=1}^{N_{\nu_t}} g_i(t) \eta_i(\mathbf{x}) \, ,
\end{equation*}
where $\eta_i(\mathbf{x})$ is the $i$-th eddy viscosity mode \blue{evaluated through a POD procedure} and $g_i(t)$ 
the corresponding coefficient.

Adding the turbulence terms to 
the SUP-ROM \eqref{reduced_system}, the new dynamical system can be written as follows:
\begin{equation}
\begin{cases}
    &\mathbf{M}\dot{\mathbf{a}} = \nu (\mathbf{B}+\mathbf{B_T}) \mathbf{a} - \mathbf{a}^T \mathbf{C} \mathbf{a} + \mathbf{g}^T (\mathbf{C}_{\text{T1}}+\mathbf{C}_{\text{T2}}) \mathbf{a} - \mathbf{H} \mathbf{b}+\tau \left( \sum_{k=1}^{N_{\text{BC}}}(U_{\text{BC},k}\mathbf{D}^k-\mathbf{E}^k \mathbf{a})\right) \, ,\\
    &\mathbf{P} \mathbf{a}=\mathbf{0} \, ,
    \end{cases}
    \label{sys_SUP_turb}
\end{equation}
where the new matrices appearing are defined as:
\begin{equation*}
    \begin{split}
    &(\mathbf{C}_{\text{T1}})_{ijk}=(\boldsymbol{\phi}_i, \eta_j \nabla \cdot \nabla \boldsymbol{\phi}_k)_{L^2(\Omega)} \, ,\\
    &(\mathbf{C}_{\text{T2}})_{ijk}=(\boldsymbol{\phi}_i, \nabla \cdot \eta_j (\nabla \boldsymbol{\phi}_k)^T)_{L^2(\Omega)}\, .
    \end{split}
\end{equation*}
When a PPE approach is considered, the FOM momentum and Poisson equation, according to the RANS turbulent model, are written in the following way:
\begin{equation*}
    \begin{cases}
    \dfrac{\partial \overline{\mathbf{u}}}{\partial t}+ \nabla \cdot (\overline{\mathbf{u}} \otimes \overline{\mathbf{u}})=\nabla \cdot \left[-\overline{p} \mathbf{I} +(\nu+\nu_t) \left(\nabla \overline{\mathbf{u}} + (\nabla \overline{\mathbf{u}})^T \right) \right] & \text{ in } \Omega \times [0,T]\, ,\\
    \Delta \overline{p}=-\nabla \cdot (\nabla \cdot (\overline{\mathbf{u}} \otimes \overline{\mathbf{u}})) +\nabla \cdot \left[ \nabla \cdot \left( \nu_t \left(\nabla \overline{\mathbf{u}} +(\nabla \overline{\mathbf{u}})^T \right) \right) \right] & \text{ in }\Omega \, ,\\
    + \text{ Boundary Conditions }& \text{ on } \Gamma \times [0,T]\, ,\\
    + \text{ Initial Conditions } & \text{ in } (\Omega,0)\, .
    \end{cases}
\end{equation*}
Consequently, the dynamical system \eqref{reduced_systemPPE} 
takes the following form:
\begin{equation}
    \begin{cases}
    \mathbf{M} \dot{\mathbf{a}}=\nu(\mathbf{B}+\mathbf{B_T})\mathbf{a}-\mathbf{a}^T \mathbf{C} \mathbf{a}+ \mathbf{g}^T (\mathbf{C}_{\text{T1}} +\mathbf{C}_{\text{T2}}) \mathbf{a}-\mathbf{H}\mathbf{b}+\tau \left( \sum_{k=1}^{N_{\text{BC}}}(U_{\text{BC},k}\mathbf{D}^k-\mathbf{E}^k \mathbf{a})\right) \, ,\\
    \mathbf{D}\mathbf{b}+ \mathbf{a}^T \mathbf{G} \mathbf{a} -\mathbf{g}^T(\mathbf{C}_{\text{T3
    }} +\mathbf{C}_{\text{T4
    }})\mathbf{a} - \nu \mathbf{N} \mathbf{a}- \mathbf{L}=\mathbf{0}\, ,
    \end{cases}
    \label{sys_PPE_turb}
\end{equation}
where:
\[
(\mathbf{C}_{\text{T3}})_{ijk}=(\nabla \chi_i, \eta_j \nabla \cdot \nabla \boldsymbol{\phi}_k)_{L^2(\Omega)}\, , \quad (\mathbf{C}_{\text{T4}})_{ijk}=(\nabla \chi_i, \nabla \cdot \eta_j(\nabla \boldsymbol{\phi}_k)^T)_{L^2(\Omega)}\, .
\]
In the dynamical systems defined in \eqref{sys_SUP_turb} and \eqref{sys_PPE_turb}, the number of unknowns is $N_u+N_{sup}$ (in \eqref{sys_SUP_turb}) and $N_u$ (in \eqref{sys_PPE_turb}) 
for velocity, $N_p$ 
for pressure, and $N_{\nu_t}$ 
for the eddy viscosity. However, the number of equations is $N_u+N_{sup}+N_p$ in \eqref{sys_SUP_turb} and $N_u+N_p$ in \eqref{sys_PPE_turb}. Thus, there are more unknowns than equations and the systems 
are not closed 
for both 
the supremizer and the Poisson approach. In order to close the systems, the eddy viscosity coefficients $[g_i(t)]_{i=1}^{N_{\nu_t}}$ can be computed considering the mapping $\mathbf{g}=f(\mathbf{a})$ through either \emph{interpolation} or \emph{regression} techniques.

In the first approach, the reduced eddy viscosity coefficients are interpolated with radial basis functions \cite{lazzaro2002radial, micchelli1986interpolation}. This technique was exploited in \cite{hijazi2020data}, following the POD-I approach \cite{wang2012comparative,walton2013reduced,salmoiraghi2018free}.
In the regression approach, which is the one used in this paper, the reduced eddy viscosity coefficients are computed starting from the velocity coefficients $[a_i]_{i=1}^{N_u}$ through a feed-forward fully-connected neural network \cite{ivagnes2021data}.
The neural network considered has two hidden layers of sizes $256$ and $64$, the activation function is \emph{ReLU}, and the learning rate used for training is $1e-5$. The loss function minimized in the training process is the difference between the output of the neural network ($\mathbf{g}_{\text{NN}}=\mathbf{f}_{\text{NN}}(\mathbf{a})$) and the values of the known coefficients of the eddy viscosity field ($\mathbf{g}=(g_i)_{i=1}^{N_{\nu_t}}$), which are found from the POD procedure.

\section{Hybrid data-driven ROMs}
\label{combined_dataROM}
This section introduces the hybrid 
data-driven approach. 
The new formulation 
fuses the purely data-driven and physically-based data-driven strategies presented 
in sections \ref{dataROM} and \ref{evROM}, respectively. The reduced system following the supremizer approach is expressed as follows:
\begin{equation}
\begin{cases}
    &\mathbf{M}\dot{\mathbf{a}} = \nu (\mathbf{B}+\mathbf{B_T}) \mathbf{a} - \mathbf{a}^T \mathbf{C} \mathbf{a} + \mathbf{g}^T (\mathbf{C}_{\text{T1}}+\mathbf{C}_{\text{T2}}) \mathbf{a} - \mathbf{H} \mathbf{b}+\tau \left( \sum_{k=1}^{N_{\text{BC}}}(U_{\text{BC},k}\mathbf{D}^k-\mathbf{E}^k \mathbf{a})\right) + \boldsymbol{\tau}\blue{_u}(\ba) \, ,\\
    &\mathbf{P} \mathbf{a}=\mathbf{0} \, ,
    \end{cases}
    \label{sys_SUP_turb_corr}
\end{equation}
where:
\begin{equation}
    \boldsymbol{\tau}\blue{_u}(\ba) = \tilde{A} \ba + \ba^T \tilde{B} \ba \, .
\end{equation}
Considering the pressure Poisson approach, the reduced system is written in the following way:
\begin{equation}
    \begin{cases}
    \mathbf{M} \dot{\mathbf{a}}=\nu(\mathbf{B}+\mathbf{B_T})\mathbf{a}-\mathbf{a}^T \mathbf{C} \mathbf{a}+ \mathbf{g}^T (\mathbf{C}_{\text{T1}} +\mathbf{C}_{\text{T2}}) \mathbf{a}-\mathbf{H}\mathbf{b}+\tau \left( \sum_{k=1}^{N_{\text{BC}}}(U_{\text{BC},k}\mathbf{D}^k-\mathbf{E}^k \mathbf{a})\right) + \boldsymbol{\tau}\blue{_u}(\ba, \bb) \, ,\\
    \mathbf{D}\mathbf{b}+ \mathbf{a}^T \mathbf{G} \mathbf{a} -\mathbf{g}^T(\mathbf{C}_{\text{T3
    }} +\mathbf{C}_{\text{T4
    }})\mathbf{a} - \nu \mathbf{N} \mathbf{a}- \mathbf{L} + \boldsymbol{\tau}\blue{_p}(\ba, \bb)=\mathbf{0}\, ,
    \end{cases}
    \label{sys_PPE_turb_corr}
\end{equation}
where, as described in section \ref{datappe}, the correction terms are: 
\begin{equation}
    \boldsymbol{\tau}_{\text{tot}} = \left(\boldsymbol{\tau}\blue{_u}(\ba,\bb), \boldsymbol{\tau}\blue{_p}(\ba,\bb) \right) = \tilde{J}_A \mathbf{ab} + \mathbf{ab}^T \tilde{J}_B \mathbf{ab}\,.
\end{equation}

\green{
\remark{
We highlight that there is an interplay between the two data-driven techniques used in the hybrid model:
\begin{itemize}
    \item the coefficient vector $\mathbf{a}(t)$ is computed at time $t$ by solving a system in which both purely and physically data-driven closure terms appear. For example, in the hybrid data-driven PPE-ROM model:
    \begin{equation}
    \begin{cases}
    \dot{\mathbf{a}}(t)=\mathbf{f}(\mathbf{a}(t), \mathbf{b}(t), \mathbf{g}(t))+\boldsymbol{\tau_u}(\mathbf{a}(t),\mathbf{b}(t))\\
    \mathbf{h}(\mathbf{a}(t), \mathbf{b}(t), \mathbf{g}(t))+\boldsymbol{{\tau}_p}(\mathbf{a}(t),\mathbf{b}(t))=0.
    \end{cases}
    \label{sys}
    \end{equation}
    The computed coefficient vector $\mathbf{a}$ is then given as input to the neural network and used to predict the eddy viscosity reduced field ($\mathbf{g}(t+1)=f(\mathbf{a}(t))$).
    \item the computed $\mathbf{g}(t+1)$ is then used to solve \eqref{sys} at time $t+1$ and find $\mathbf{a}(t+1)$ and $\mathbf{b}(t+1)$, and so on.
\end{itemize}
}
}

\section{Numerical Results}
\label{results}

\begin{figure}[h!]
\includegraphics[width = \textwidth, trim={0 74cm 0 0 }, clip]{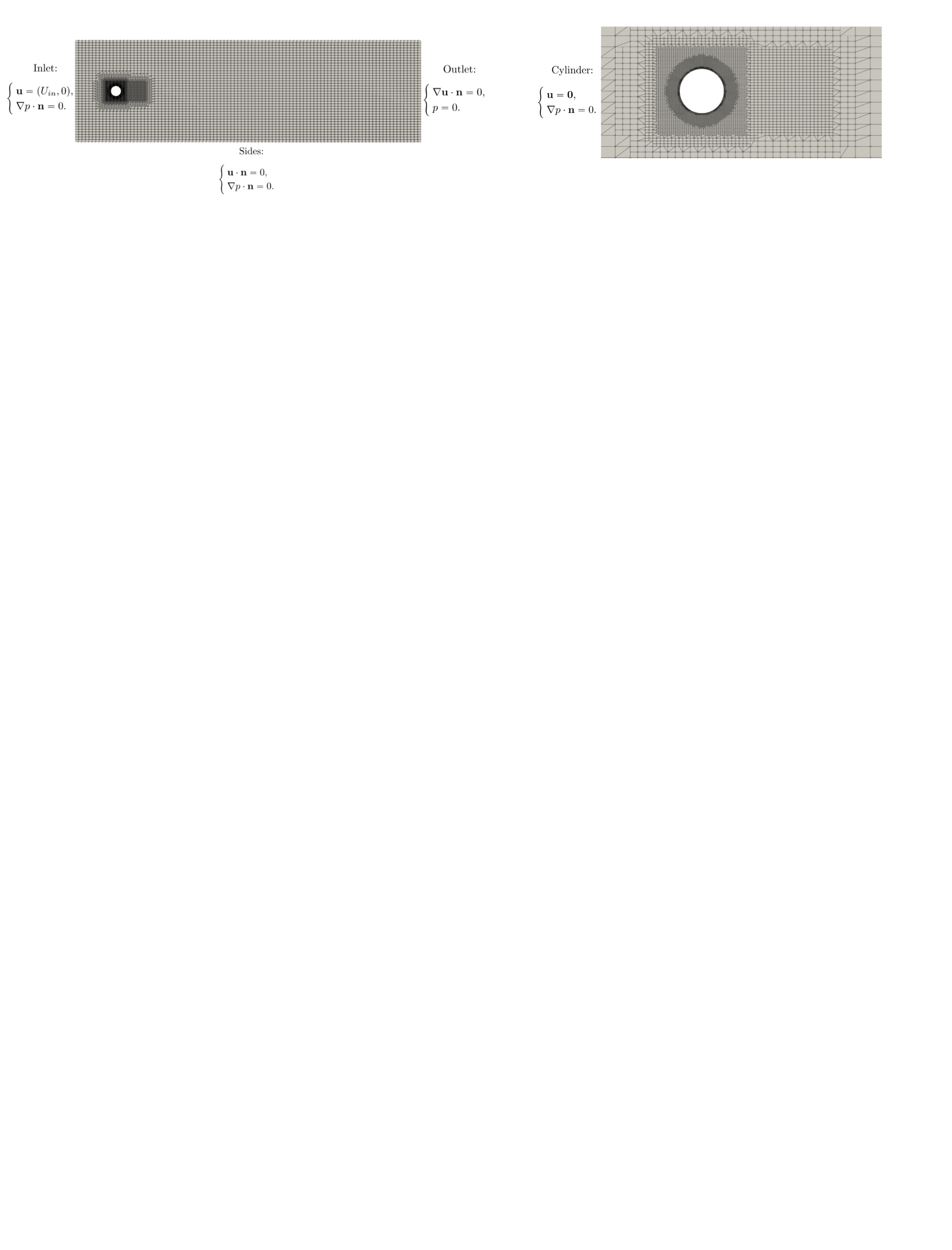}
\caption{The mesh used in simulations (left) and
the mesh zoomed around the cylinder (right), with the corresponding boundary conditions.
}
\label{mesh}
\end{figure}

The case study considered to test the data-driven 
ROMs described in sections \ref{dataROM}, \ref{evROM} and \ref{combined_dataROM} is an unsteady flow past a circular cylinder. The case is studied in two dimensions and the mesh used 
has $11644$ cells. The mesh and the boundary conditions set for velocity and pressure are 
displayed in Figure \ref{mesh} \cite{hijazi2020data}. The diameter of the cylinder is $D= \SI{1}{\metre}$, the fluid kinematic viscosity $\nu=\SI{1e-4}{\metre \squared \per \second}$, and the velocity at the inlet is horizontal and fixed: $U_{in}=\SI{5}{\metre \per \second}$.


For the offline stage, the software \emph{OpenFOAM} is used to obtain the full order fields, using the 
unsteady solver pimpleFoam and the $\kappa-\omega$ model for the turbulence treatment. The pimpleFoam solver is based on a PIMPLE approach for pressure velocity coupling with under-relaxation technique, which consists 
of the coupling of a SIMPLE \cite{patankar2018numerical} and a PISO strategy \cite{issa1986solution}.
The number of time instances considered in the offline phase is $5000$, considering one FOM snapshot every $\SI{0.004}{\second}$ and starting from the second $20$ from the beginning of the simulation. Since the time step of the full order simulation is set to $\SI{0.0002}{\second}$, we are using \emph{undersampling}, taking one snapshots every 20 time steps. 

Then, the 
POD is performed. The POD modes for the velocity, pressure, and 
supremizer fields are extracted from the numerical data, making use of the library \emph{ITHACA-FV} \cite{ithacasite, stabile2018finite, Stabile2017CAIM}. 

The simulations for the online ROM procedure are computed in a dedicated \emph{Python} script, and both the SUP-ROM and the PPE-ROM approaches are tested.

The results of the high-fidelity simulations are compared to the results obtained by solving the reduced order dynamical systems with or without the data-driven terms and considering the parameter $\tau=1000$, which was introduced in system \eqref{reduced_system}. 
More specifically, the results are evaluated in terms of the percentage errors of the reduced pressure and velocity fields with respect to the correspondent high-fidelity fields, evaluated at each time step, in the $L^2$ norm.
Since the projection of the full order solution on the reduced POD space is the best possible result which can be achieved with a given 
number of modes, the solution of the reduced system cannot improve with respect to that projection. 
For this reason, 
the percentage errors with respect to the projection of the full order solution on the reduced POD space are used in the numerical investigation.
The percentage errors with respect to the full order fields at each $j$-th time step are evaluated in the following way:
\begin{equation}
    \varepsilon_u(t_j)=\dfrac{||\mathbf{u}_r^{abs}(\mathbf{x}, t_j)-\mathbf{u}_d^{abs} (\mathbf{x}, t_j)||_{L^2(\Omega)}}{||\mathbf{u}^{abs}_d (\mathbf{x}, t_j)||_{L^2(\Omega)}}\, , \quad \varepsilon_p(t_j)=\dfrac{||p_r(\mathbf{x}, t_j)-p_d(\mathbf{x}, t_j)||_{L^2(\Omega)}}{||p_d(\mathbf{x}, t_j)||_{L^2(\Omega)}}\, ,
    \label{errors1}
\end{equation}
where \blue{the following quantities are considered}:
\begin{itemize}
   \blue{ \item the reduced velocity and pressure fields:
    $$\mathbf{u}_r(\mathbf{x}, t_j)=\sum_{i=1}^{r} a_i(t_j) \boldsymbol{\phi}_i(\mathbf{x}) , \quad p_r(\mathbf{x},t_j)=\sum_{i=1}^{q} b_i(t_j) \chi_i(\mathbf{x})  ,$$
    where the coefficients $a_i(t_j)$ and $b_i(t_j)$ are the solutions of the dynamical systems \eqref{supgen} (in the supremizer approach) and \eqref{ppegen} (in the Poisson approach);
    \item the approximated full order fields of velocity
    and pressure, which are evaluated starting from the first $d$ modes, where $d=100$ when a supremizer approach is considered, and $d=50$ when a pressure Poisson approach is considered. Such a term reads 
    $$\mathbf{u}_d(\mathbf{x}, t_j)=\sum_{i=1}^{d} a^{\text{snap}}_i(t_j) \boldsymbol{\phi}_i (\mathbf{x}) , \quad p_d(\mathbf{x}, t_j)=\sum_{i=1}^{d_p} b^{\text{snap}}_i (t_j) \chi_i(\mathbf{x})  .$$
    We remark that when the supremizer approach is considered,  $\{\boldsymbol{\phi}_i\}_{i=51}^{100}=\{\mathbf{s}_i(\chi_i)\}_{i=1}^{50}$ are the supremizer modes.}

\end{itemize}

Considering the supremizer approach, the most general form for the dynamical system solved at each time step is the following:
\begin{equation}
\begin{cases}
    \mathbf{M} \dot{\mathbf{a}}^{i}=\nu(\mathbf{B}+\mathbf{B_T})\mathbf{a}^i-(\mathbf{a}^{i})^T \mathbf{C} \mathbf{a}^i-\mathbf{H}\mathbf{b}^i+  c_{u} \boldsymbol{\tau}_{u}(\mathbf{a}^i, \mathbf{b}^i)&+  \\ +
     c_{t} \left( (\mathbf{g}^i)^T (\mathbf{C}_{T1}+\mathbf{C}_{T2})\mathbf{a}^i \right) \blue{+\tau \left( \sum_{k=1}^{N_{\text{BC}}}(U_{\text{BC},k}\mathbf{D}^k-\mathbf{E}^k \mathbf{a}^i)\right)} &\text{ at each }i=1,...,M\, ,\\
    \mathbf{P}\mathbf{a}^i =\mathbf{0} &\text{ at each }i=1,...,M\, ,
    \end{cases}
    \label{supgen}
\end{equation}
where $M$ is the total number of time steps in the online phase. 

Considering the pressure Poisson approach, the reduced dynamical system can be written as follows:
\begin{equation}
\begin{cases}
    \mathbf{M} \dot{\mathbf{a}}^{i}=\nu(\mathbf{B}+\mathbf{B_T})\mathbf{a}^i-(\mathbf{a}^{i})^T \mathbf{C} \mathbf{a}^i-\mathbf{H}\mathbf{b}^i + c_u \boldsymbol{\tau}^u(\mathbf{a}^i, \mathbf{b}^i)  &+\\
    +c_{t}  \left( (\mathbf{g}^i)^T (\mathbf{C}_{T1}+\mathbf{C}_{T2})\mathbf{a}^i \right) \blue{+\tau \left( \sum_{k=1}^{N_{\text{BC}}}(U_{\text{BC},k}\mathbf{D}^k-\mathbf{E}^k \mathbf{a}^i)\right)}&\text{ at each }i=1,...,M\, ,\\
   \mathbf{D}\mathbf{b}^i+ (\mathbf{a}^i)^T \mathbf{G} \mathbf{a}^i - \nu \mathbf{N} \mathbf{a}^i- \mathbf{L}+ c_{t}  \left( (\mathbf{g}^i)^T (\mathbf{C}_{T3}+\mathbf{C}_{T4})\mathbf{a}^i \right)& +\\ +c_p \boldsymbol{\tau}^{p}(\mathbf{a}^i, \mathbf{b}^i)=\mathbf{0} &\text{ at each }i=1,...,M \, .
    \end{cases}
    \label{ppegen}
\end{equation}

The matrices and tensors appearing in the two formulations are specified in sections \ref{sec:g-rom}, \ref{dataROM}, and \ref{evROM}.
In systems \eqref{supgen} and \eqref{ppegen}, the parameters $c_u$, $c_p$, and $c_t$ are introduced to either include or combine the different data-driven 
strategies in the reduced formulation. \medskip

For the time integration, two different approaches are used. 
Specifically, the time derivative $\dot{\mathbf{a}}^i$ appearing in the momentum equation in \eqref{supgen} and \eqref{ppegen} is computed following (i) an Euler first order time scheme, and (ii) an implicit second order time scheme.
It is worth remarking that the second order time scheme corresponds to the scheme implemented in \emph{OpenFOAM} and used to solve the full order problem.
\medskip

In this section, the numerical results are analyzed 
with respect to three different criteria:
\begin{itemize}
    \item in the first part (section \ref{error_integrals}), the accuracy of the data-driven ROMs is analyzed
    by varying the number of ROM modes 
    between $1$ and $10$;
    \item the second and the third sections (\ref{sup_turb_res} and \ref{ppe_turb_res}) 
    investigate the accuracy of the data-driven ROMs by first fixing the number of modes to $5$ for the velocity, pressure, and supremizer fields, and then varying the time. 
    Sections \ref{sup_turb_res} and \ref{ppe_turb_res} investigate the supremizer and the Poisson approaches, respectively;
    \item the last section (\ref{graph_sec}) is dedicated to the graphical comparison of the 
    velocity and pressure fields produced by the data-driven ROMs.
\end{itemize}

The following cases are analyzed and compared, for both the first and the second order time schemes, and for both SUP-ROM and PPE-ROM approaches:
\begin{itemize}
    \item $c_t=c_u=0$ (and $c_p=0$ for PPE-ROM), i.e., standard ROM without the addition of any data-driven term;
    \item $c_t=0$, $c_u=1$ (and $c_p=1$ in PPE-ROM), i.e., purely data-driven model;
    \item $c_t=1$, $c_u=0$ (and $c_p=0$ in PPE-ROM), i.e., 
    physically-based data-driven model;
    \item $c_t=c_u=1$ (and $c_p=1$ in PPE-ROM), i.e., hybrid data-driven 
    model.
\end{itemize}

In sections \ref{sup_turb_res}, \ref{ppe_turb_res}, and \ref{graph_sec}, the results provided by the combination of the two different data-driven techniques are compared to those obtained 
in the previous works \cite{paper1} and \cite{hijazi2020data}, where only one of these  two techniques is used. 

\subsection{Hybrid data-driven approach in different modal regimes}
\label{error_integrals}
In this section, a comparison 
of different modal regimes is carried out for both the supremizer and the pressure Poisson approaches. 
The models considered for the investigation are:
\begin{itemize}
    \item the SUP-ROM, with and without the closure turbulence model and the extra velocity correction term;
    \item the PPE-ROM, with and without the 
    turbulence model, 
    extra velocity correction term (in the momentum equation), and 
    pressure correction term (in the pressure Poisson equation).
\end{itemize}
The accuracy is evaluated in terms of \blue{percentage integral errors} for the absolute value of velocity and pressure on a time window of $500$ time instants. The expressions for the overall time window $L^2$ errors are the following:
\blue{
\begin{equation*}
    \dfrac{\int_0^T||\mathbf{u}^{abs}_r(\mathbf{x}, t_j)-\mathbf{u}^{abs}_d (\mathbf{x}, t_j)||_{L^2(\Omega)} \, dt }{ \int_0^T ||\mathbf{u}^{abs}_d (\mathbf{x}, t_j)||_{L^2(\Omega)} \, dt} \times 100 \, , \quad \dfrac{\int_0^T||p_r(\mathbf{x}, t_j)-p_d(\mathbf{x}, t_j)||_{L^2(\Omega)} \, dt}{\int_0^T||p_d(\mathbf{x}, t_j)||_{L^2(\Omega)} \, dt} \times 100\,
\end{equation*}}
where $T$ corresponds to $\SI{2}{\second}$.
The errors defined above are analyzed for all the 
cases for both the first and the second order time schemes.

\begin{figure}[h!]
\centering
\subfloat[Percentage error of velocity ($1^{st}$ order)]{\includegraphics[ height=5.3cm]{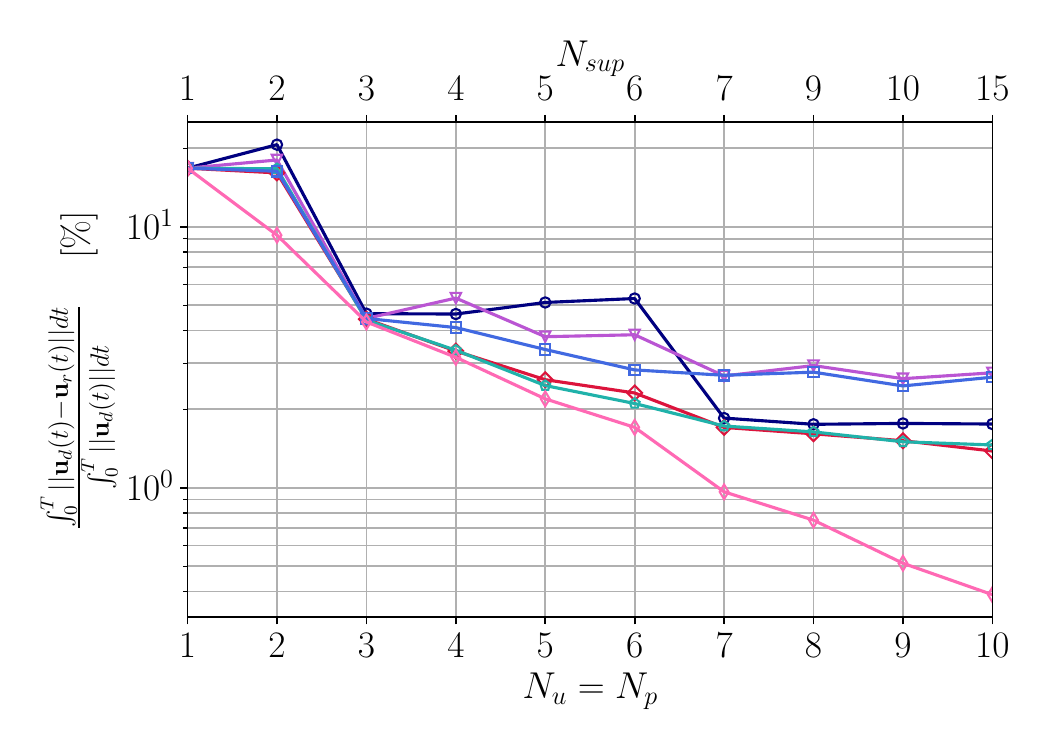} \label{1U}}
\subfloat[Percentage error of pressure ($1^{st}$ order)]{\includegraphics[ height=5cm ]{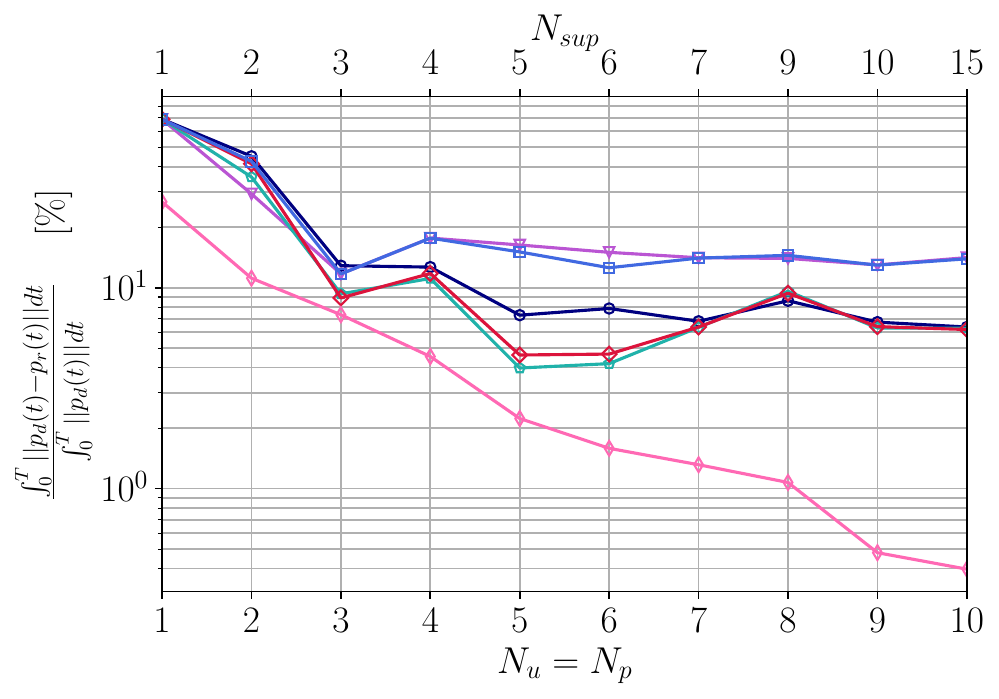} \label{1p}}\\
\subfloat[Percentage error of velocity ($2^{nd}$ order)]{\includegraphics[ height=5.3cm]{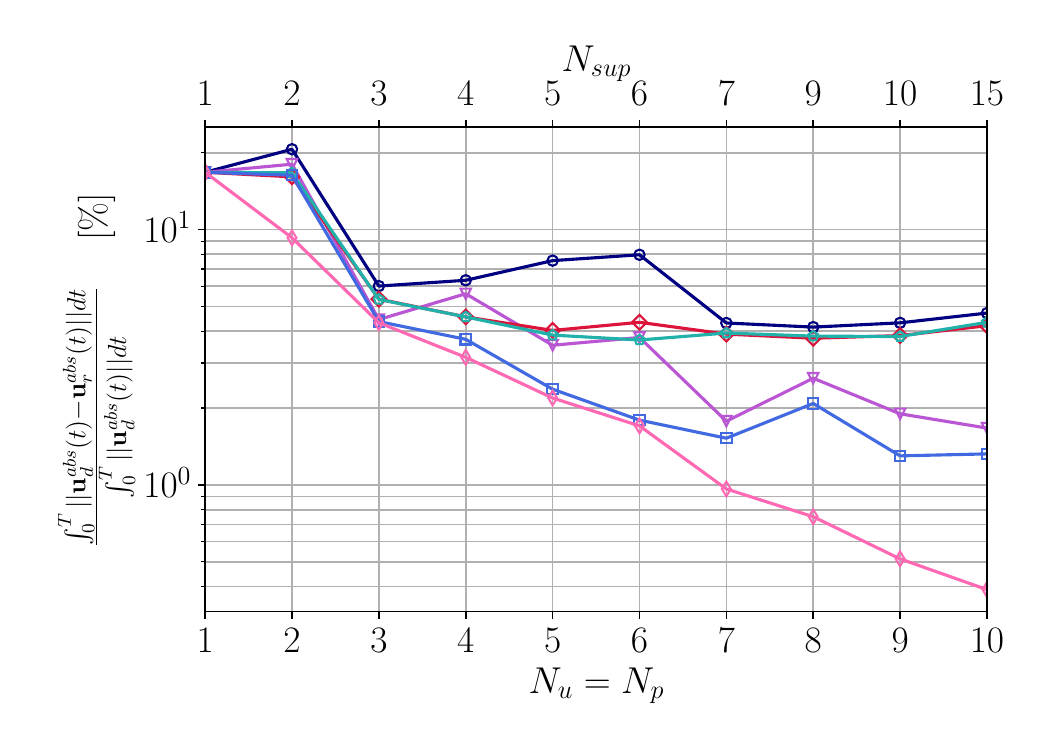} \label{2U}}
\subfloat[Percentage error of pressure  ($2^{nd}$ order)]{\includegraphics[ height=5cm]{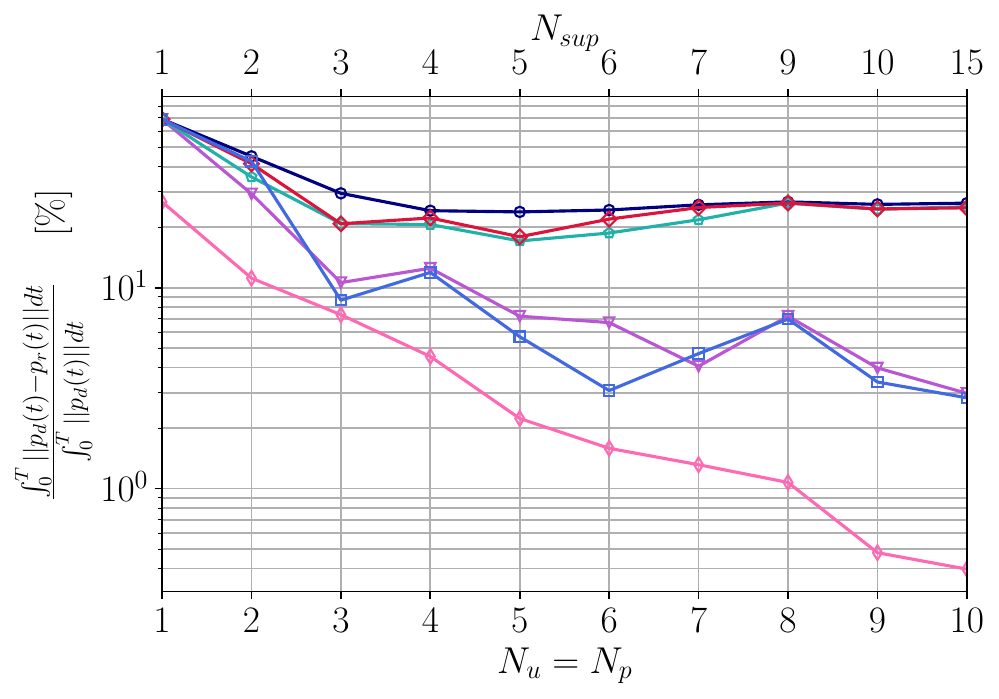} \label{2p}}
\caption{\blue{Percentage integral errors} of the absolute value of velocity and pressure, 
varying the number of modes. The model is the SUP-ROM with a first (\protect\subref{1U}, \protect\subref{1p}) and second (\protect\subref{2U}, \protect\subref{2p}) order time scheme. The cases represented are the following: without any data-driven term (\includestandalone{markers-plots/mark0}); physically-based data-driven model (\includestandalone{markers-plots/mark1}); purely data-driven model for the velocity (\includestandalone{markers-plots/mark2} and \includestandalone{markers-plots/mark3} unconstrained and constrained, respectively); hybrid model (\includestandalone{markers-plots/mark4}); projection (\includestandalone{markers-plots/mark5}).
}
\label{modes_SUP}
\end{figure}

\begin{figure}[h!]
\centering
\subfloat[Percentage error of velocity ($1^{st}$ order)]{\includegraphics[ height=5.3cm]{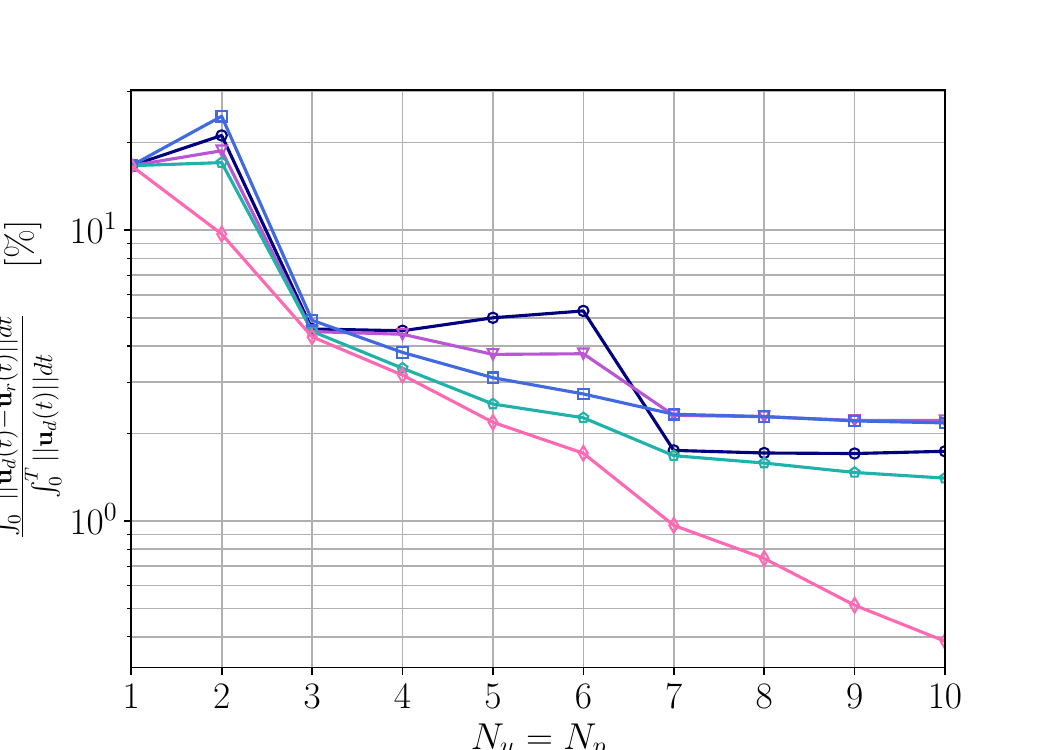} \label{1U_p}}
\subfloat[Percentage error of pressure ($1^{st}$ order)]{\includegraphics[ height=5cm ]{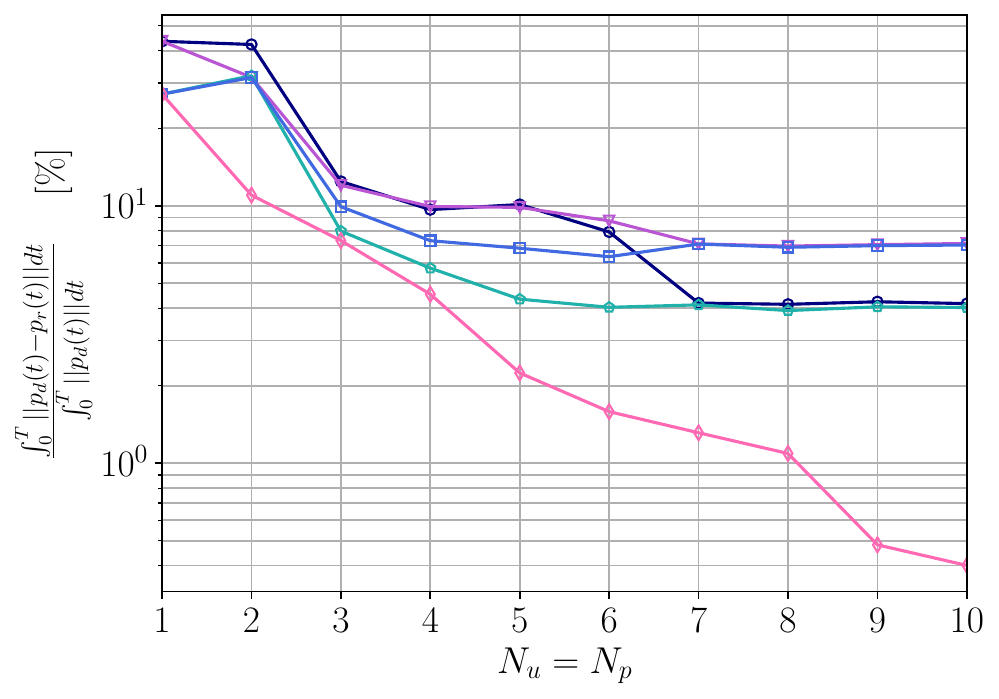} \label{1p_p}}\\
\subfloat[Percentage error of velocity ($2^{nd}$ order)]{\includegraphics[ height=5.3cm]{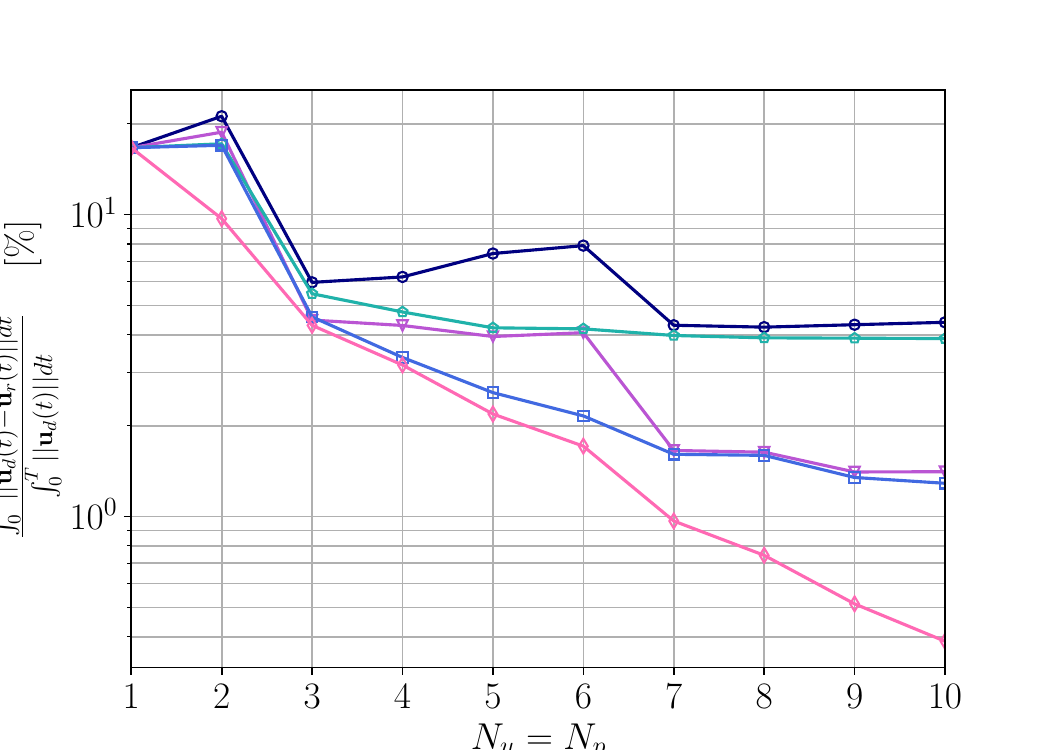} \label{2U_p}}
\subfloat[Percentage error of pressure  ($2^{nd}$ order)]{\includegraphics[ height=5cm]{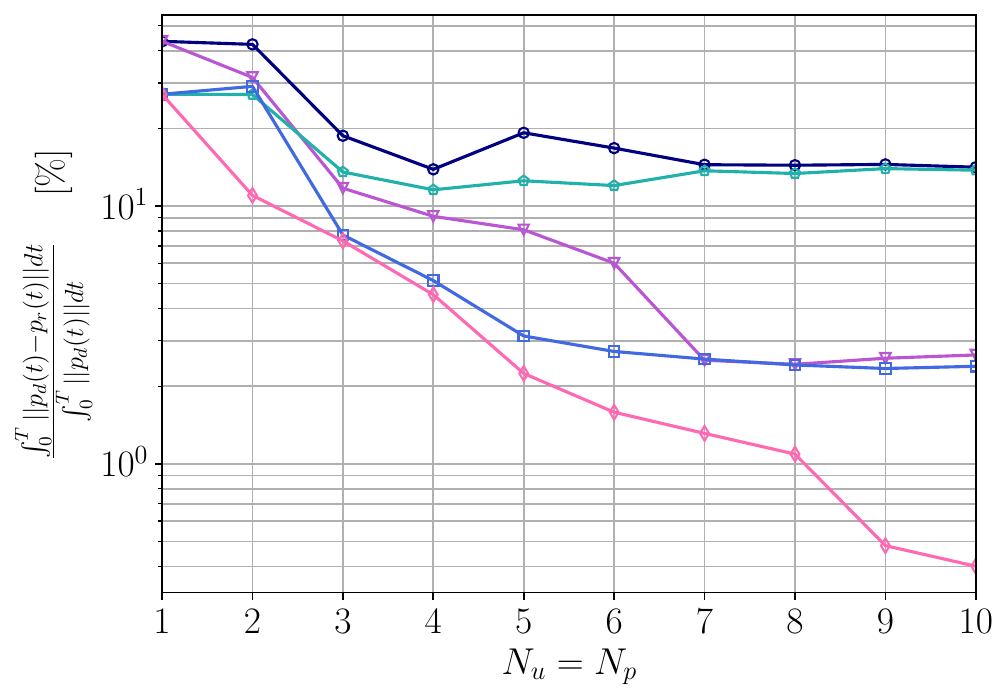} \label{2p_p}}
\caption{\blue{Percentage integral errors} of the absolute value of velocity and of pressure, varying the number of modes.
The model is the PPE-ROM with a first (\protect\subref{1U_p}, \protect\subref{1p_p}) and second (\protect\subref{2U_p}, \protect\subref{2p_p}) order time scheme. The cases represented are the following: without any data-driven term (\includestandalone{markers-plots/mark0}); physically-based data-driven model (\includestandalone{markers-plots/mark1}); purely data-driven model (\includestandalone{markers-plots/mark3}); hybrid data-driven model (\includestandalone{markers-plots/mark4}); projection (\includestandalone{markers-plots/mark5}).}
\label{modes_PPE}
\end{figure}

The graphical results are showed in semi-logarithmic plots in Figure \ref{modes_SUP} for the supremizer approach, and in Figure \ref{modes_PPE} for the Poisson approach.

In Figure \ref{modes_SUP}, the number of modes considered 
in the supremizer approach 
satisfies $N_{sup}>N_p$ when $N_u=N_p$ is equal to $8,\,9$, and $10$, in order to avoid stability issues. From Figures \ref{modes_SUP} and \ref{modes_PPE}, one can note that the data-driven corrections have a more significant positive effect when a small number of modes is considered; the improvement with respect to the standard Galerkin-ROM, especially for the pressure field, is not visible when the number of modes is larger than $8$.

In both the supremizer and the Poisson approaches, the 
time integration scheme plays an important role 
in the evaluation of the 
data-driven ROM accuracy.
In fact, when a first order time scheme is used, the introduction of the turbulence treatment in the formulation does not improve the results with respect to the standard ROM and 
the results obtained in 
\cite{paper1} with the DD-VMS-ROM. This fact is particularly evident for a number of modes between $7$ and $10$ (Figures \ref{modes_SUP}\protect\subref{1U},\subref{1p} and \ref{modes_PPE}\protect\subref{1U_p},\subref{1p_p}).

\blue{When the standard ROM is considered, the results obtained with a first-order time scheme outperform the results obtained with a second-order time scheme: this is likely due to the reduced numerical dissipation associated with the second-order time scheme, which makes the system more exposed to numerical instabilities.

When the turbulence modeling is added to the reduced formulation, results improve only if we consider a second-order time scheme, i.e., if we have consistency with respect to the full order model.}

The analysis 
in the following section focuses on $N_u=N_p=5$ (and $N_{sup}=5$ for the supremizer enrichment), i.e., in the marginally-resolved regime.

 \subsection{Comparison of data-driven VMS-SUP-ROMs}
 \label{sup_turb_res}
 
In this section, the results of different data-driven techniques are compared for the supremizer approach for a first and a second order time schemes. 

 \begin{figure}[h!]
\centering
\subfloat[Percentage error of velocity]{\includegraphics[ scale=0.45]{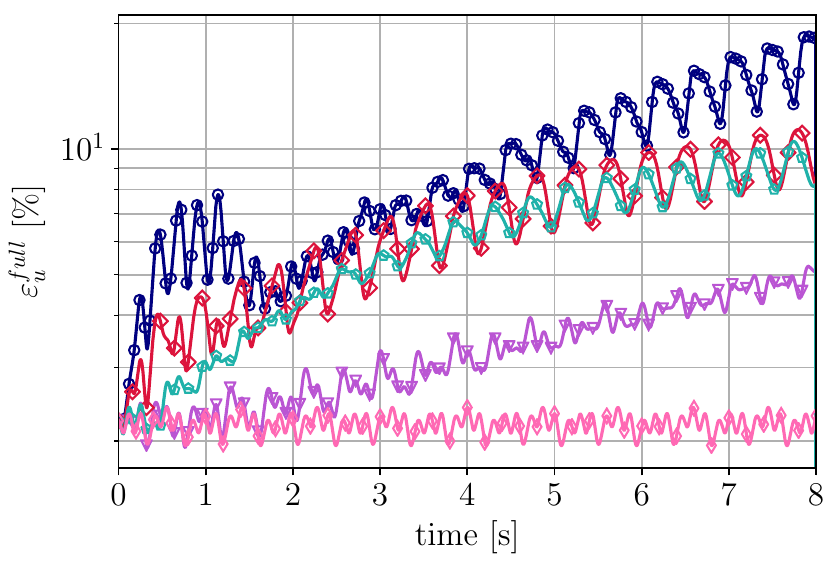} \label{sup_t_1U}}
\subfloat[Percentage error of pressure]{\includegraphics[ scale=0.45]{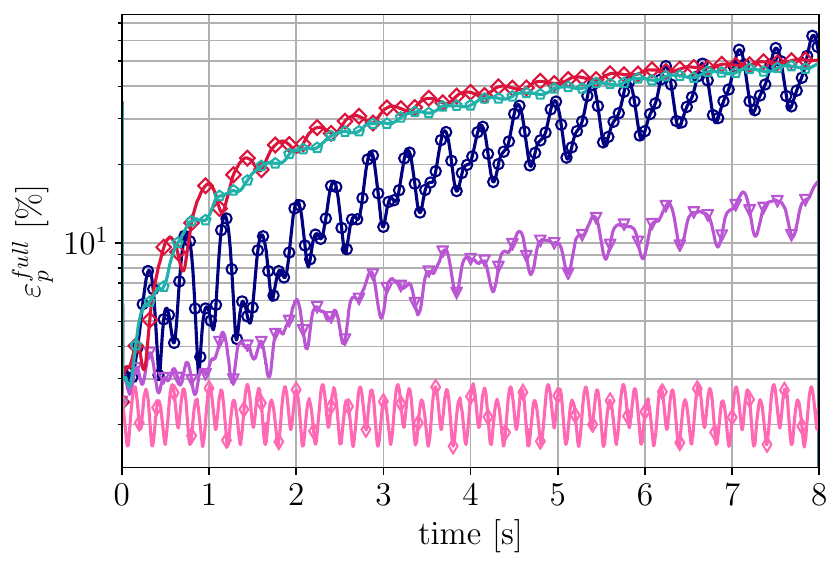} \label{sup_t_1p}}\\
\subfloat[Percentage error of velocity]{\includegraphics[ scale=0.45]{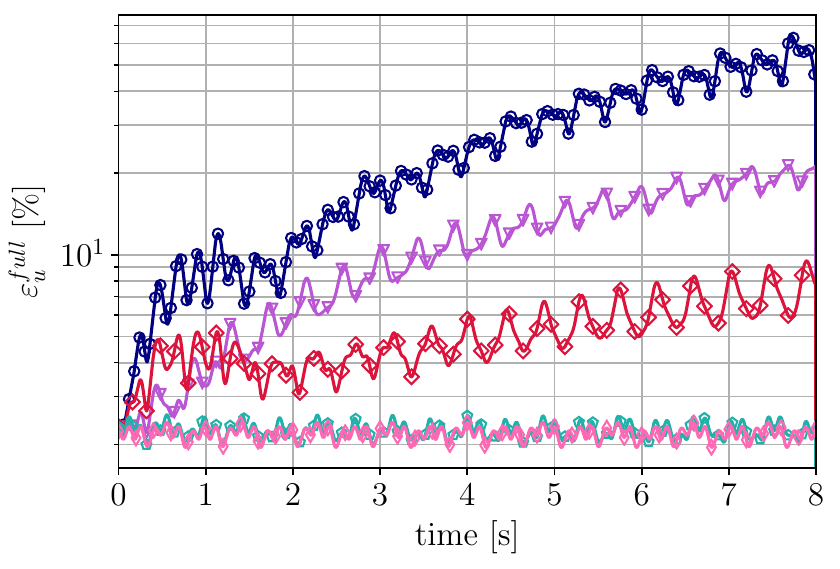} \label{sup_t_2U}}
\subfloat[Percentage error of pressure]{\includegraphics[ scale=0.45]{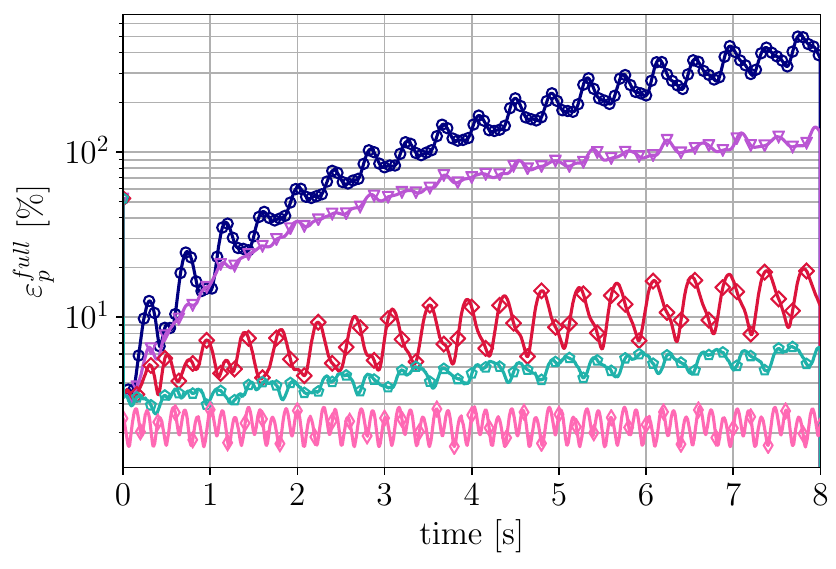} \label{sup_t_2p}}
\caption{Percentage errors of the absolute value of velocity and pressure, considering $N_u=N_p=N_{sup}=5$. The model is the SUP-ROM with a first (\protect\subref{sup_t_1U}, \protect\subref{sup_t_1p}) and second (\protect\subref{sup_t_2U}, \protect\subref{sup_t_2p}) order time scheme. Results include the following cases: without any data-driven term (\includestandalone{markers-plots/mark0}); purely data-driven constrained model for the velocity (\includestandalone{markers-plots/mark1}); physically-based data-driven model (\includestandalone{markers-plots/mark2}); hybrid data-driven model (\includestandalone{markers-plots/mark3}); projection (\includestandalone{markers-plots/mark5}).}
\label{5modes_sup}
\end{figure}

For the simulations of the reduced systems \eqref{sys_SUP_turb} and \eqref{sys_SUP_turb_corr}, the coefficients of the reduced eddy viscosity field $(g_i)_{i=1}^{N_{\nu_t}}$ are computed 
by using a fully-connected neural network, starting from the velocity coefficients $(a_i)_{i=1}^{N_u}$. The network is composed of 
two hidden layers, the ReLU function is used as an activation function in the network, and the learning rate is set to $10^{-5}$.

 The momentum equation correction term 
 is obtained by solving the constrained optimization problem \eqref{opt_problem2}, since it provides the best performance 
 with respect to the velocity accuracy for a low number of modes, as pointed out in section \ref{datasup}. 
 
 Figures \ref{5modes_sup} \protect\subref{sup_t_1U} and \protect\subref{sup_t_1p} display the results obtained using a first order time scheme. In 
 this case, the inclusion of a turbulence model does not appear to have a completely positive effect on accuracy.  
 In addition, coupling the turbulence and correction 
 strategies does not lead to higher accuracy, especially when the pressure field is considered. In such a case in fact, the hybrid method leads to worse results than in the no-correction case, as already pointed out in section \ref{error_integrals}.

When considering a second order time 
discretization scheme (Figures \ref{5modes_sup} \protect\subref{sup_t_2U} and \protect\subref{sup_t_2p}), the results obtained with turbulence modeling or with both correction closure terms and turbulence terms  are more accurate than the results of the standard ROM, and are 
close to the 
projection of the full order solution on the reduced POD space.
In particular, 
in terms of accuracy of the velocity reduced field, the results are very similar to the projected field.

\subsection{Comparison of data-driven VMS-PPE-ROMs}
\label{ppe_turb_res}
\begin{figure}[h!]
\centering
\subfloat[Percentage error of velocity]{\includegraphics[ scale=0.45]{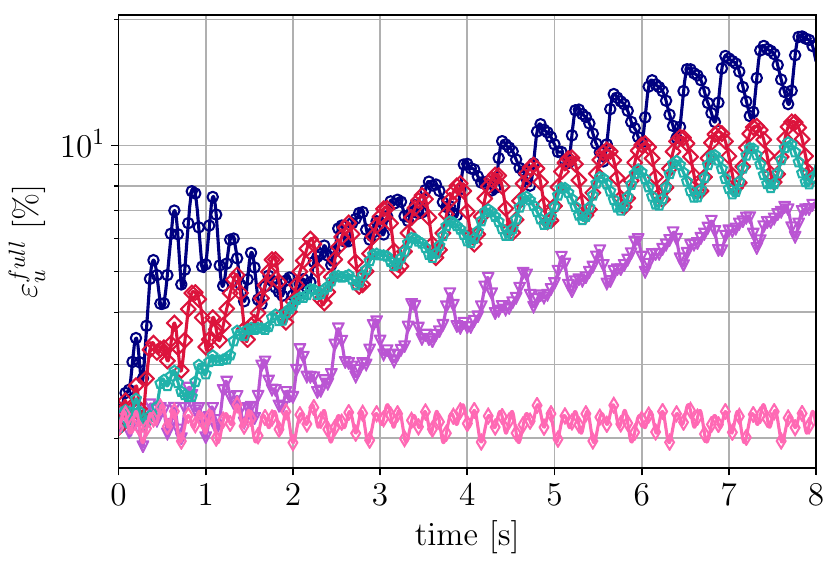} \label{ppe_t_1U}}
\subfloat[Percentage error of pressure]{\includegraphics[ scale=0.45]{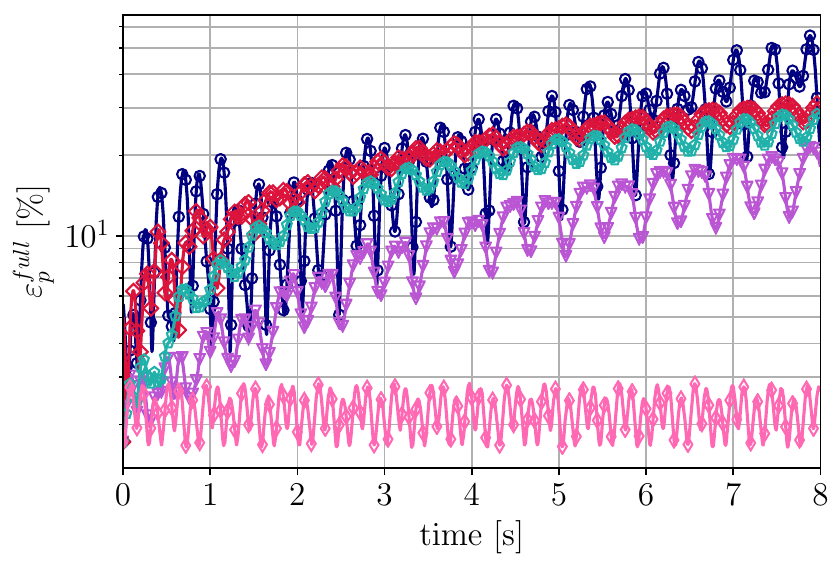} \label{ppe_t_1p}}\\
\subfloat[Percentage error of velocity]{\includegraphics[ scale=0.45]{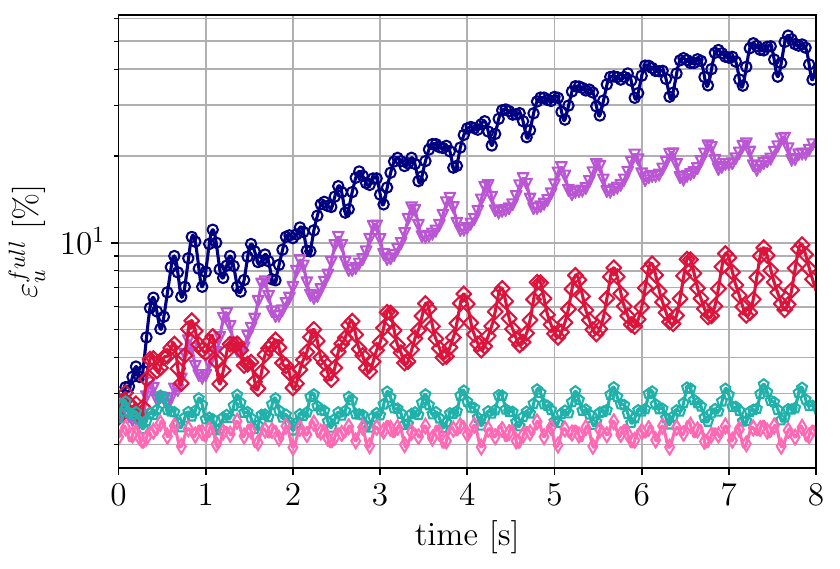} \label{ppe_t_2U}}
\subfloat[Percentage error of pressure]{\includegraphics[ scale=0.45]{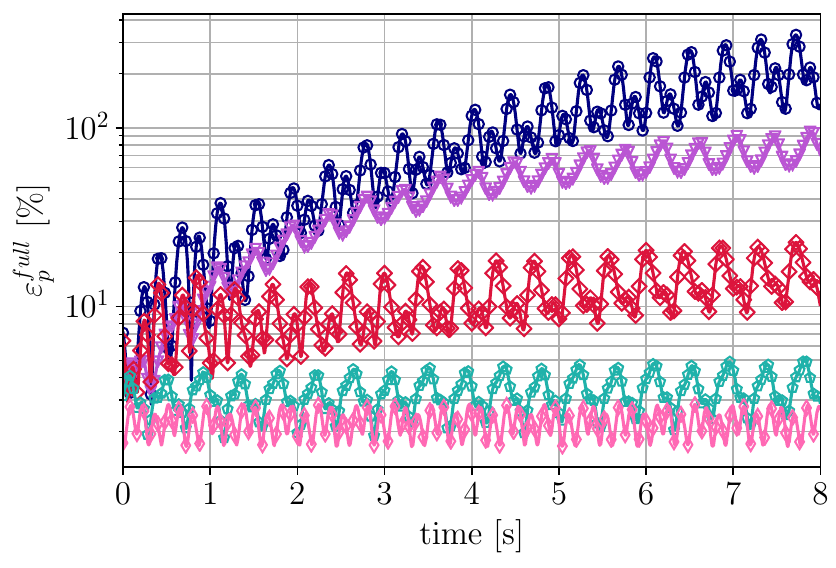} \label{ppe_t_2p}}
\caption{Percentage errors of the absolute value of velocity and pressure, considering $N_u=N_p=N_{sup}=5$. The model is the PPE-ROM, with a first (\protect\subref{ppe_t_1U}, \protect\subref{ppe_t_1p}) and second (\protect\subref{ppe_t_2U}, \protect\subref{ppe_t_2p}) order time scheme. Results include the following cases: without any data-driven term (\includestandalone{markers-plots/mark0}); purely data-driven model, for both velocity and pressure (\includestandalone{markers-plots/mark1}); physically-based data-driven model (\includestandalone{markers-plots/mark2}); hybrid data-driven model (\includestandalone{markers-plots/mark3}); projection (\includestandalone{markers-plots/mark5}).}
\label{5modes_ppe}
\end{figure}

In this section, the combined effect of data-driven terms and turbulence modeling is 
evaluated for the PPE approach for a simulation lasting $8$ seconds in Figure \ref{5modes_ppe}.
The correction terms are built starting from the first $2$ seconds and all corrections are 
constructed by using the approach detailed in Section \ref{datappe}.

The differences between 
Figures \ref{5modes_ppe} \protect\subref{ppe_t_1U},\protect\subref{ppe_t_1p} and \ref{5modes_ppe} \protect\subref{ppe_t_2U},\protect\subref{ppe_t_2p} confirm what observed 
in the SUP-ROM investigation. 
When a first order scheme is used (Figure \ref{5modes_ppe}\protect\subref{ppe_t_1U},\protect\subref{ppe_t_1p}), the hybrid 
data-driven approach does not significantly improve 
the accuracy of the eddy viscosity or purely data-driven approaches.
However, when a second order time scheme is 
used (Figure \ref{5modes_ppe} \protect\subref{ppe_t_2U},\protect\subref{ppe_t_2p}), the 
combination of turbulence modeling and corrections leads to 
a significant accuracy improvement.
\blue{\remark{Different behaviours of the reduced system in Figures \ref{modes_SUP}, \ref{modes_PPE}, \ref{5modes_sup} and \ref{5modes_ppe} when using a first or second order time scheme, can be justified with the different dissipation associated with the two numerical time integration schemes.

In fact, from the numerical examples, we deduce that, when considering a ROM which includes all the stabilizations of the FOM (i.e. data-driven + eddy viscosity), it is important to use the same temporal scheme employed at the FOM level (i.e. second order scheme). In the case of a ROM with missing stabilization terms (i.e. missing eddy viscosity for example) a first order scheme introduces additional numerical stabilization and gives a beneficial contribution.}}

We also note that the model has an excellent extrapolation efficiency. 
Indeed, although the corrections are constructed with data from the first $2$ seconds, they increase the ROM accuracy on the interval $[2,8]$ seconds.
Finally, we point out that the instability of the second order integration scheme is damped by the addition of the turbulence model.

Overall, the numerical results in this section show that adding the data-driven velocity and pressure corrections proposed in \cite{paper1} can significantly increase the accuracy of the turbulence ROMs proposed in \cite{hijazi2020data}.

\subsection{Qualitative results}
\label{graph_sec}
The inclusion of 
correction terms and 
turbulence modeling in the reduced formulations is also examined from a graphical point of view for the SUP-ROM and PPE-ROM approaches. 
 The results are graphically represented on the test case grid 
 by using the open-source application \emph{Paraview} and the results are compared to those obtained in \cite{paper1} and \cite{hijazi2020data}.

The second order time integration scheme is used since it provides the best results in Sections \ref{sup_turb_res} and \ref{ppe_turb_res}.

The POD is performed on \blue{the snapshots evaluated at the time instances within} the interval $[79.992,99.992]$ seconds and the reduced order systems \eqref{supgen} and \eqref{ppegen} are solved in the interval $[79.992,87.992]$ seconds, since the maximum length of the online simulations carried out is $8$ seconds. For this reason, all the fields are captured at the final time step of online simulations, 
i.e., at $87.992$ seconds.

\begin{figure}[h!]
\centering
    \subfloat[]{\includegraphics[width=0.46\textwidth]{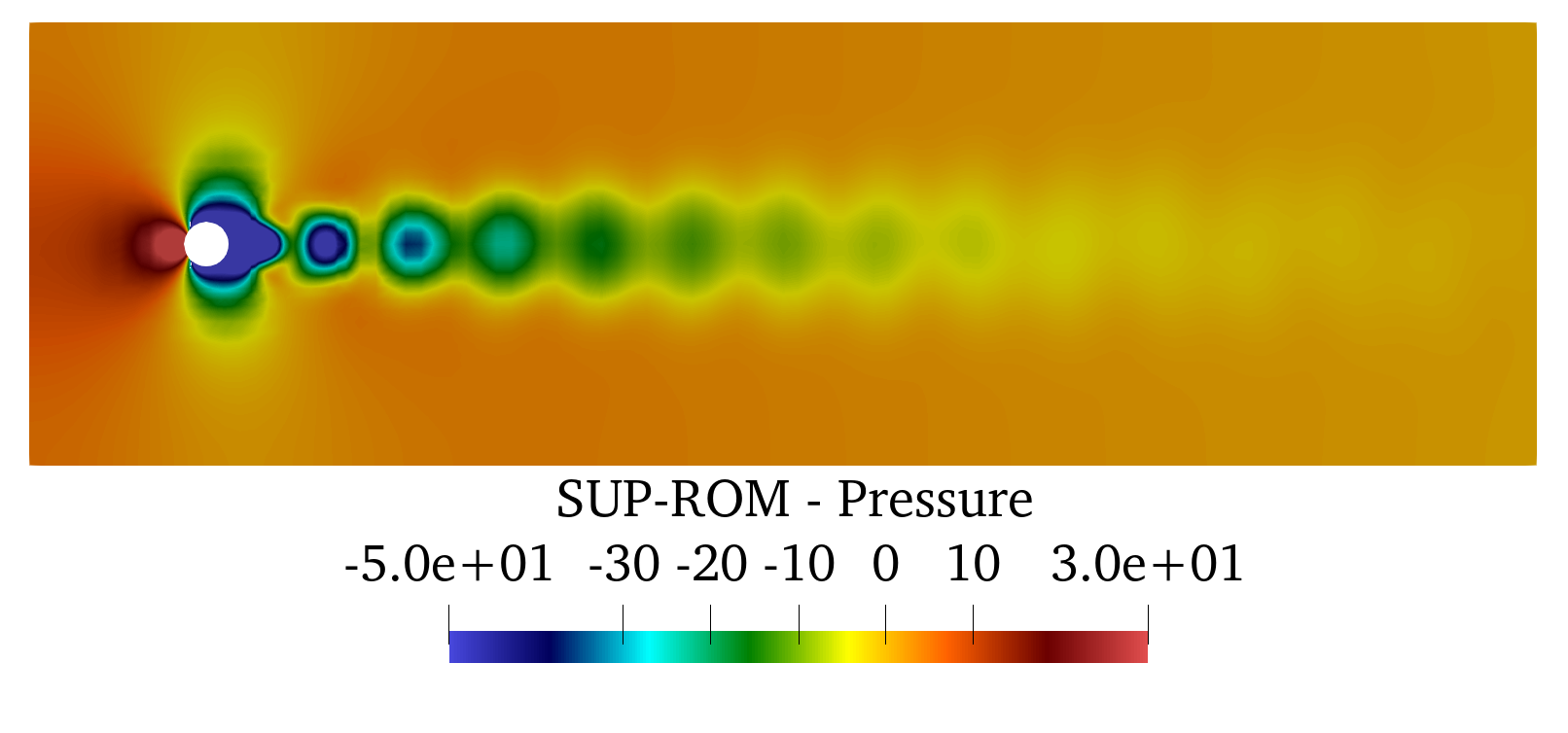}}
    \subfloat[]{\includegraphics[width=0.46\textwidth]{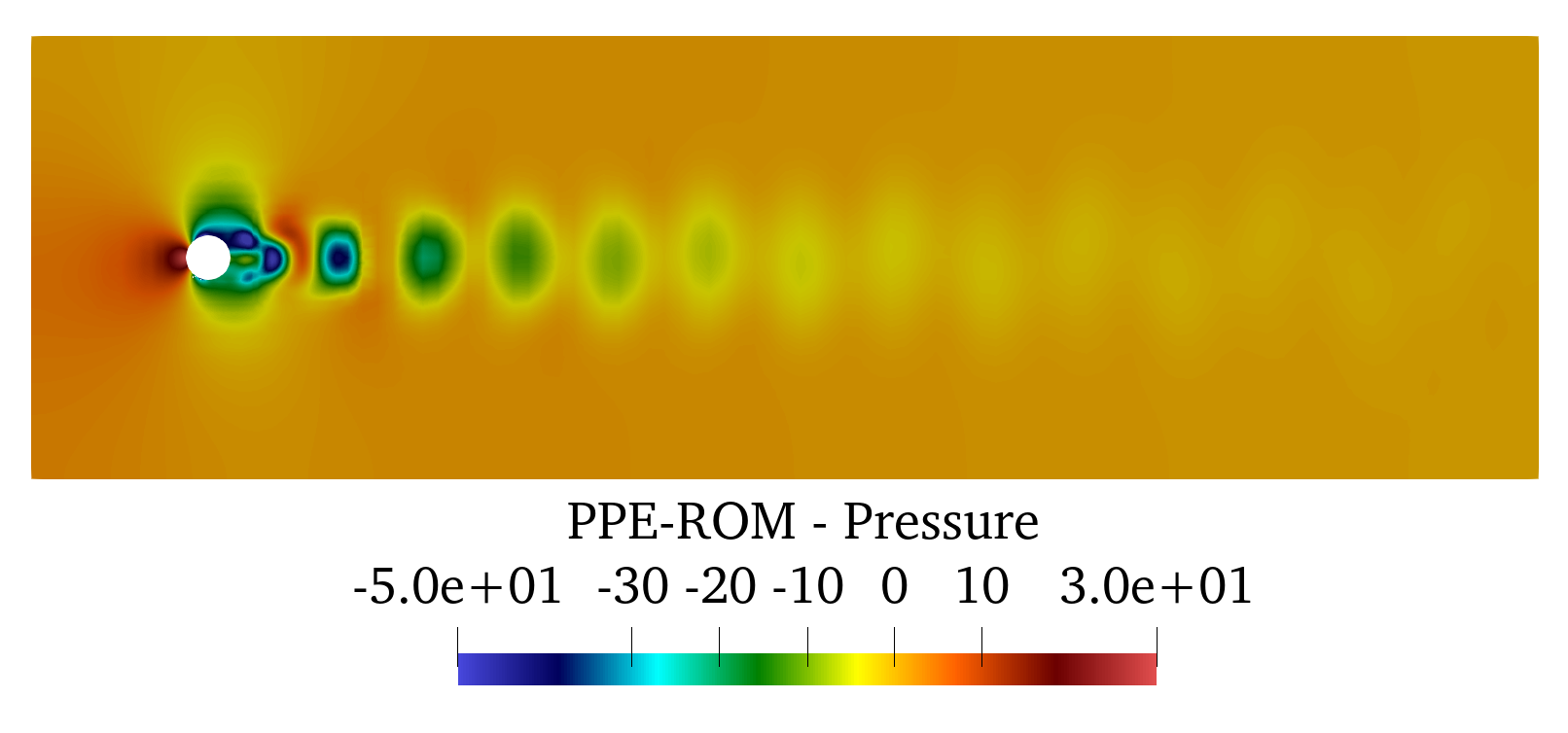}}
    \\ \vspace{-0.3cm}
    \subfloat[]{\includegraphics[width=0.46\textwidth]{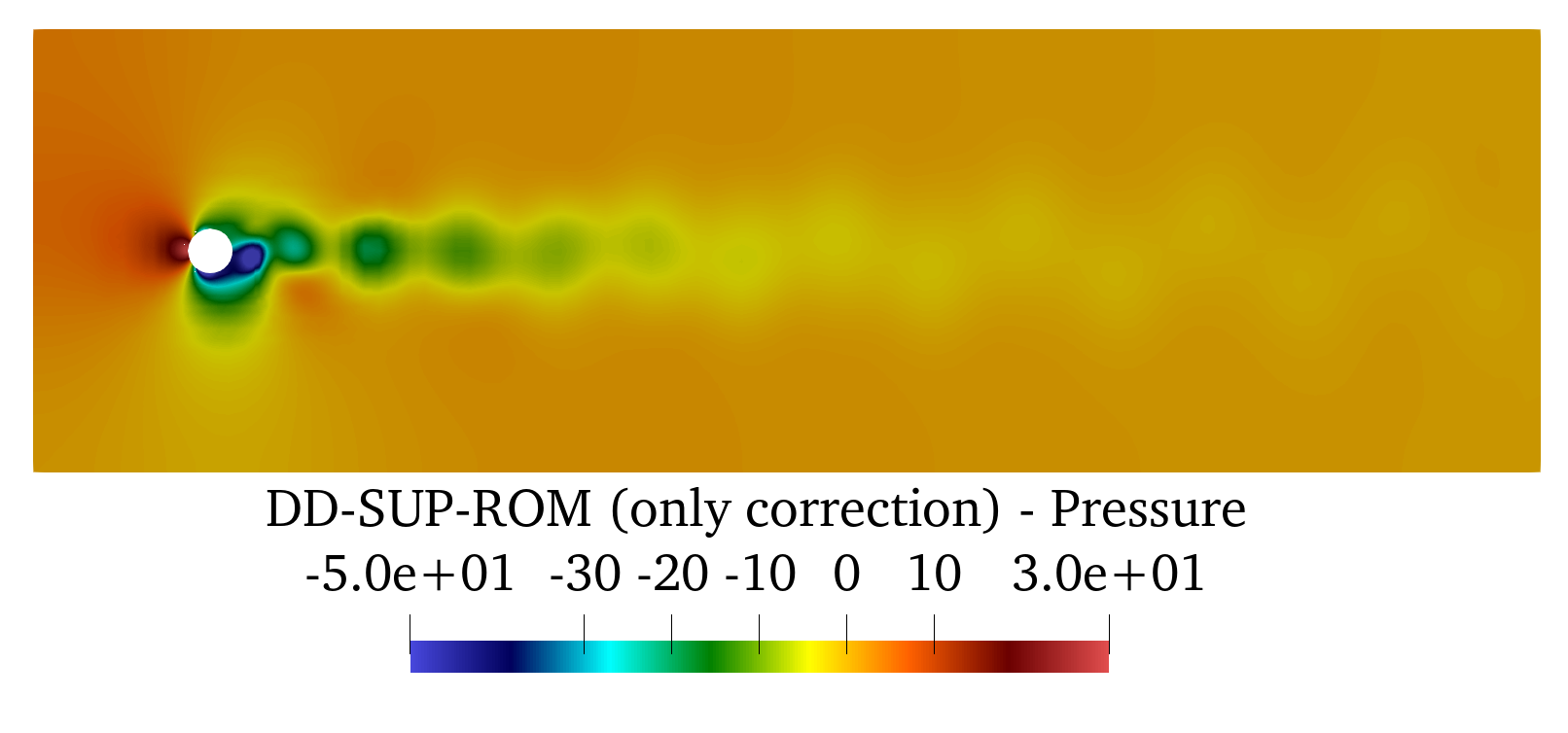}}
    \subfloat[]{\includegraphics[width=0.46\textwidth]{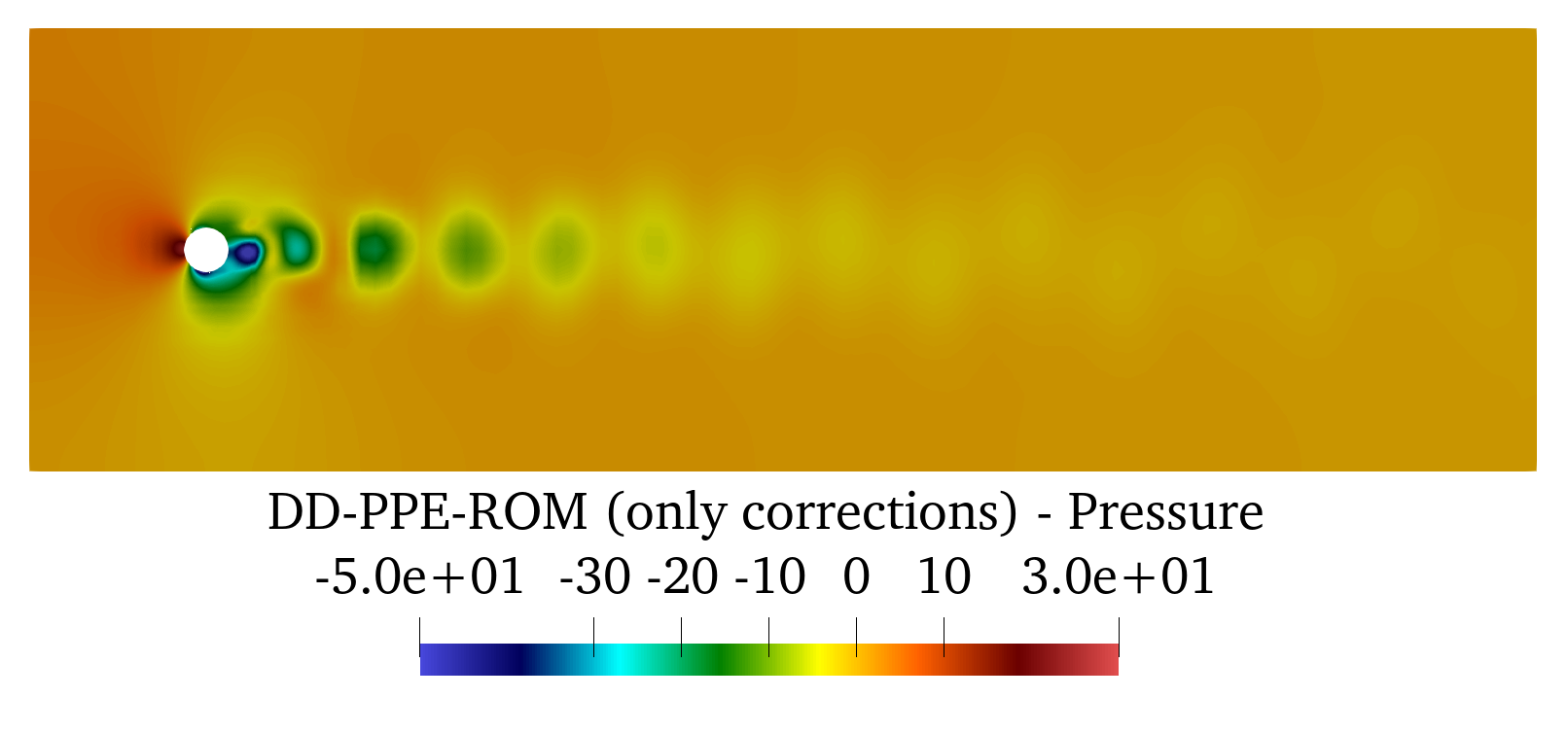}}\\ \vspace{-0.3cm}
    \subfloat[]{\includegraphics[width=0.46\textwidth]{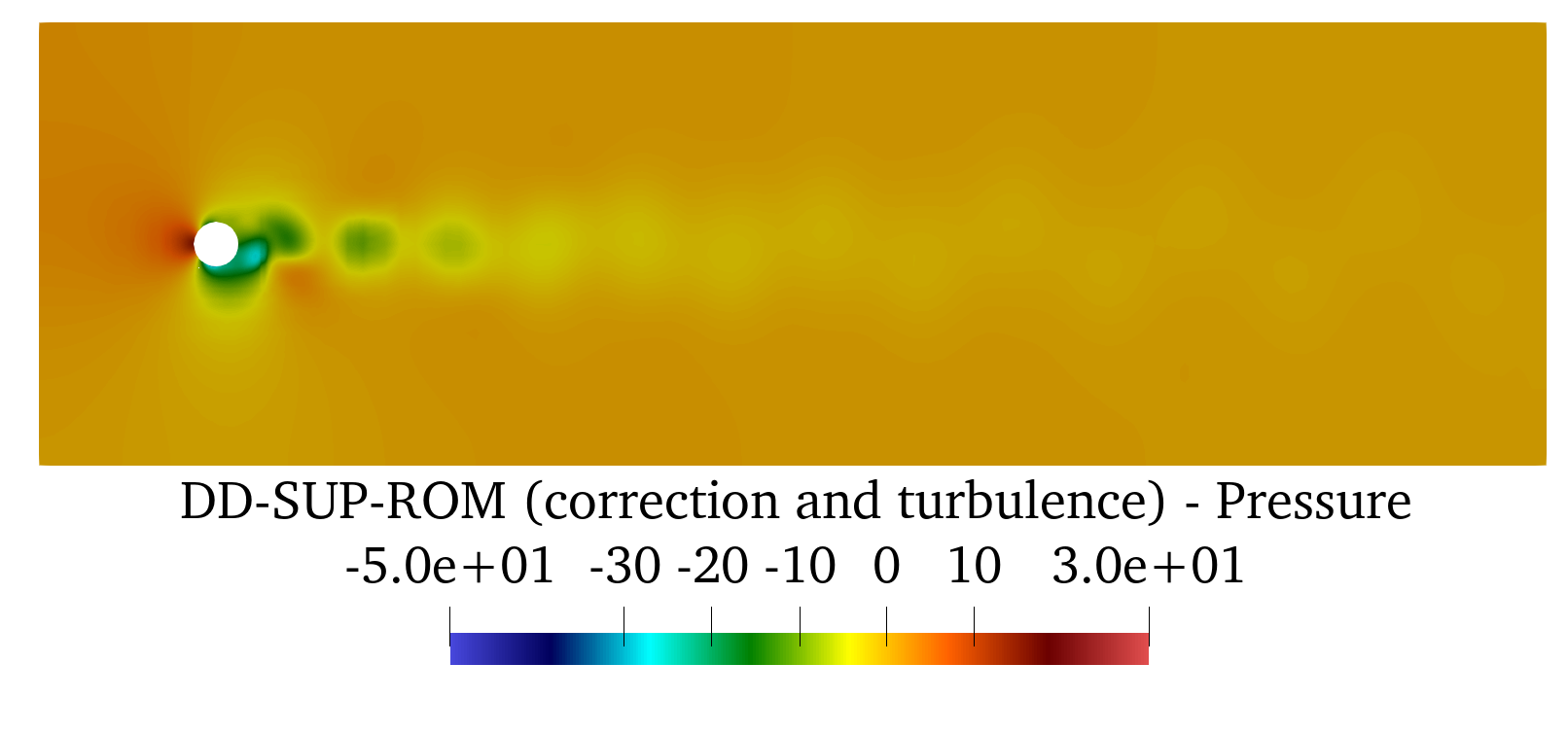}}
    \subfloat[]{\includegraphics[width=0.46\textwidth]{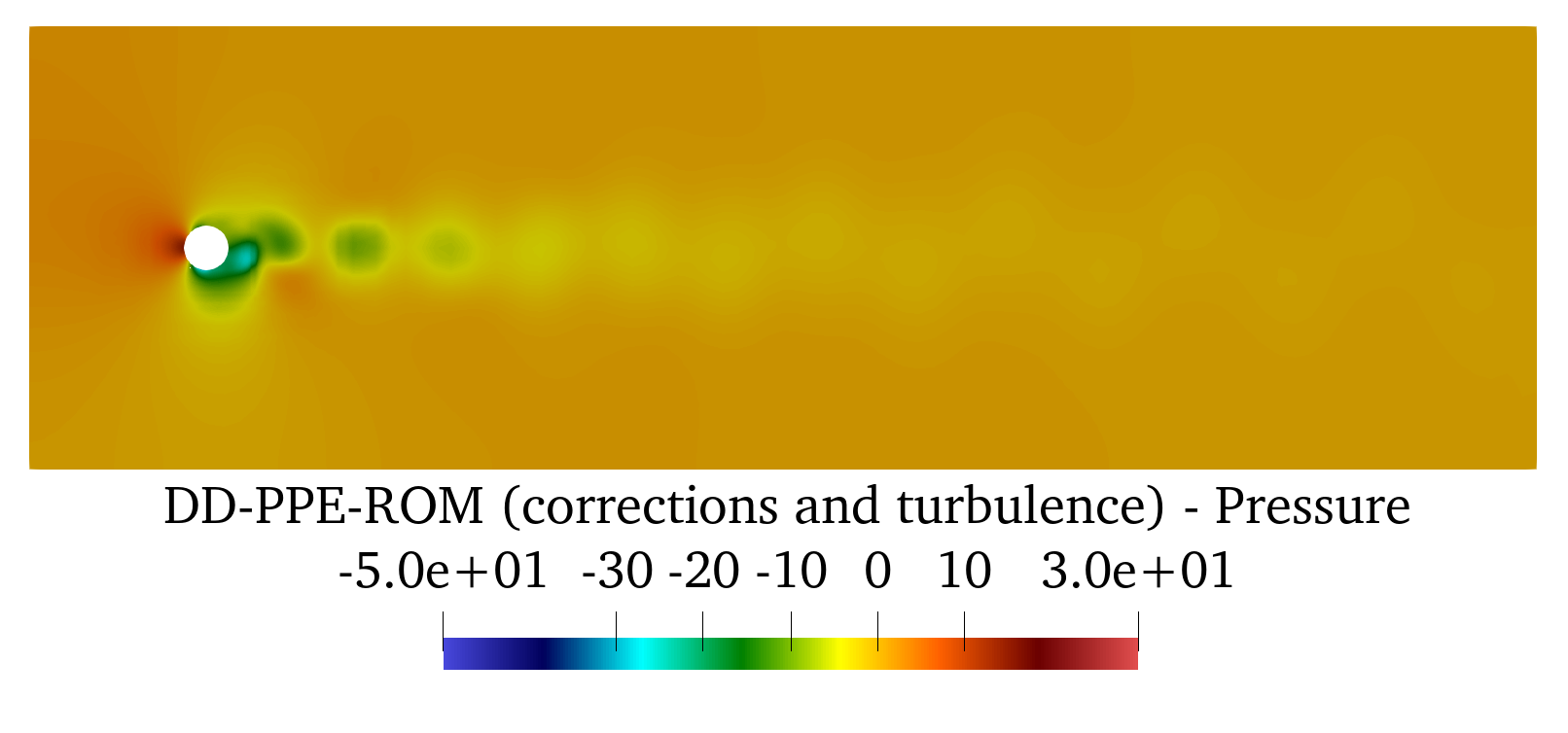}}\\ \vspace{-0.5cm}
    \subfloat[]{\includegraphics[width=0.46\textwidth]{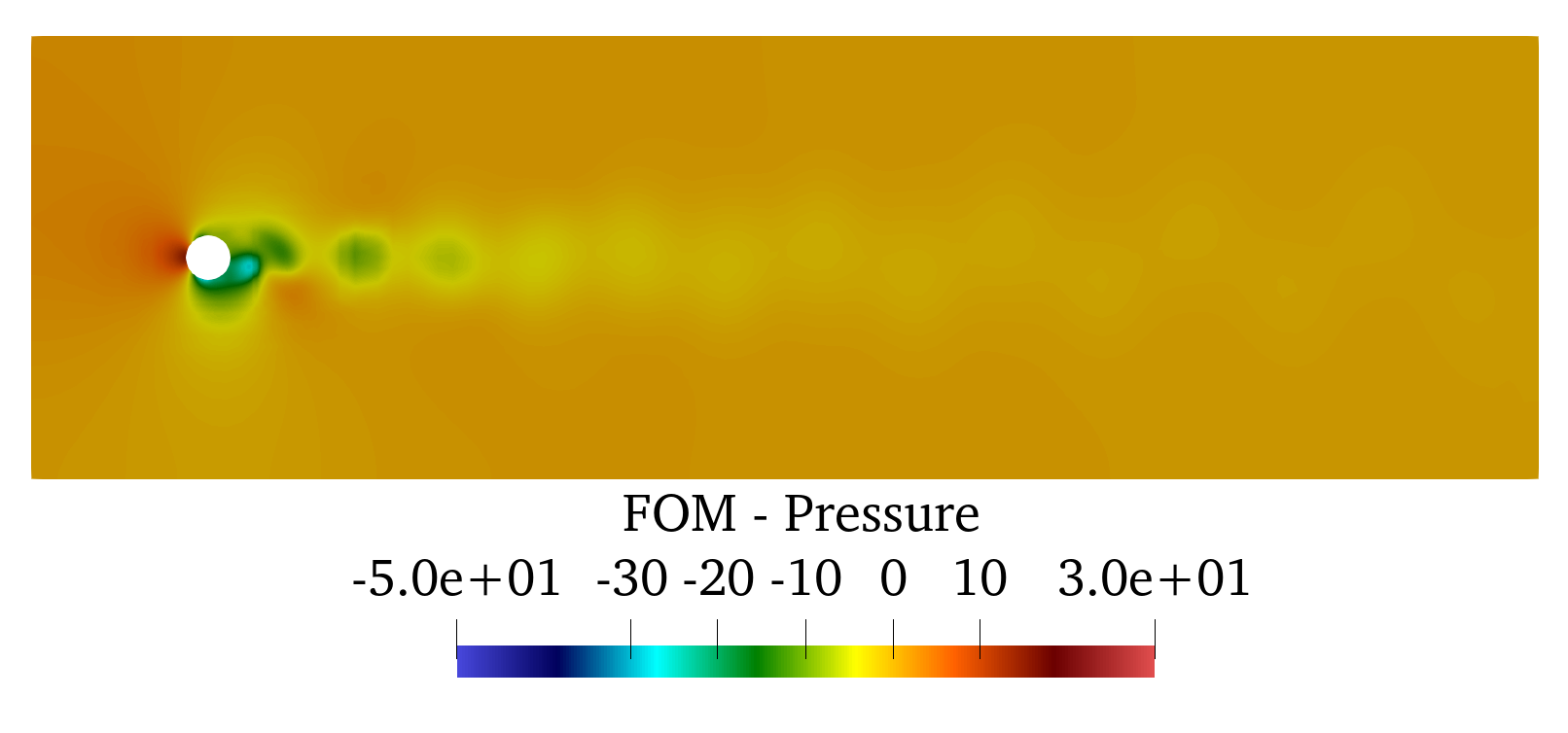}}
    \vspace{-0.4cm}
    \caption{Representation of the pressure field for the FOM, the SUP-ROM and the PPE-ROM simulations with and without the data-driven terms.}
    \label{paraview1}
\end{figure}

\begin{figure}[h!]
\centering
    \subfloat[]{\includegraphics[width=0.46\textwidth]{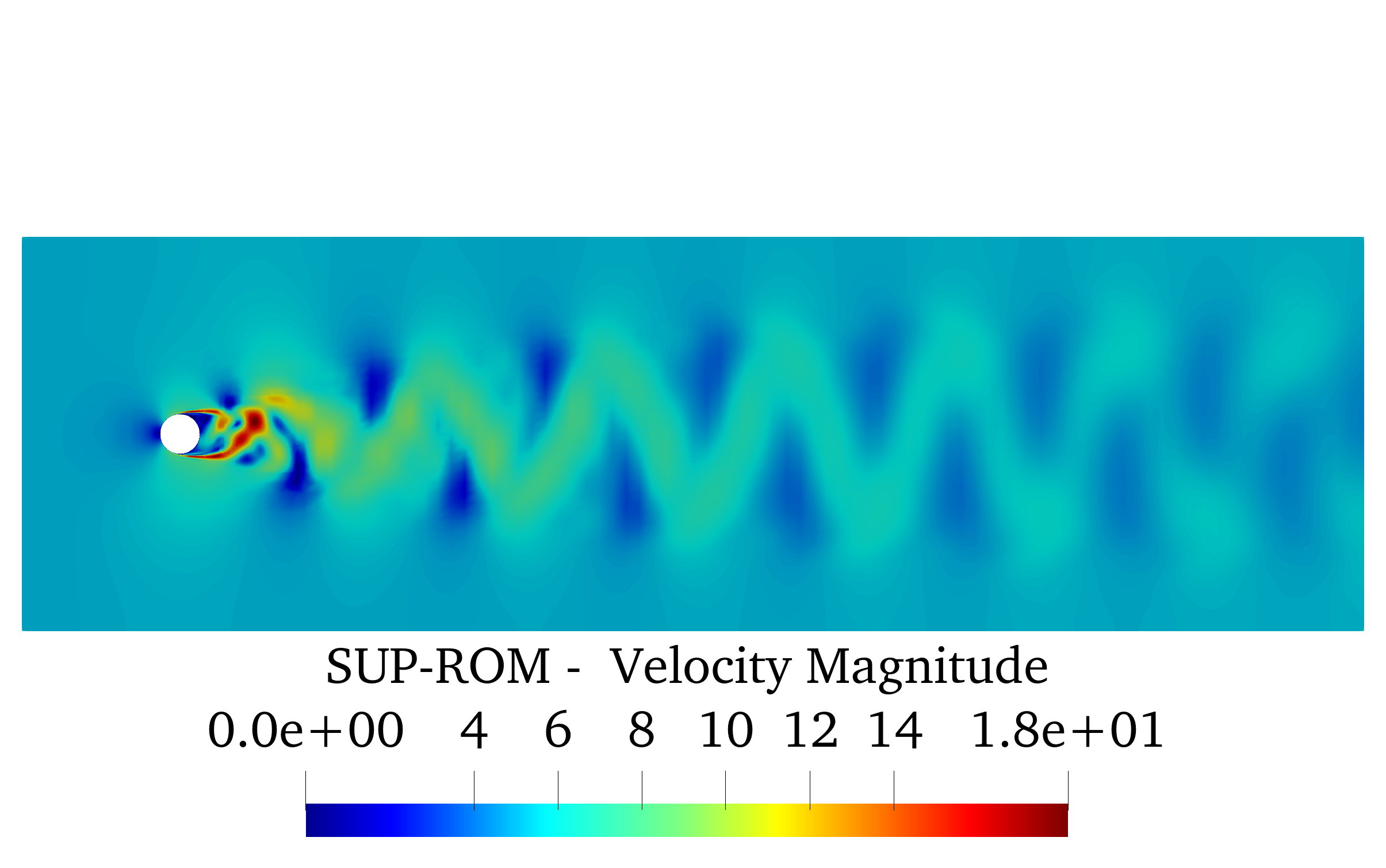}}
    \subfloat[]{\includegraphics[width=0.46\textwidth]{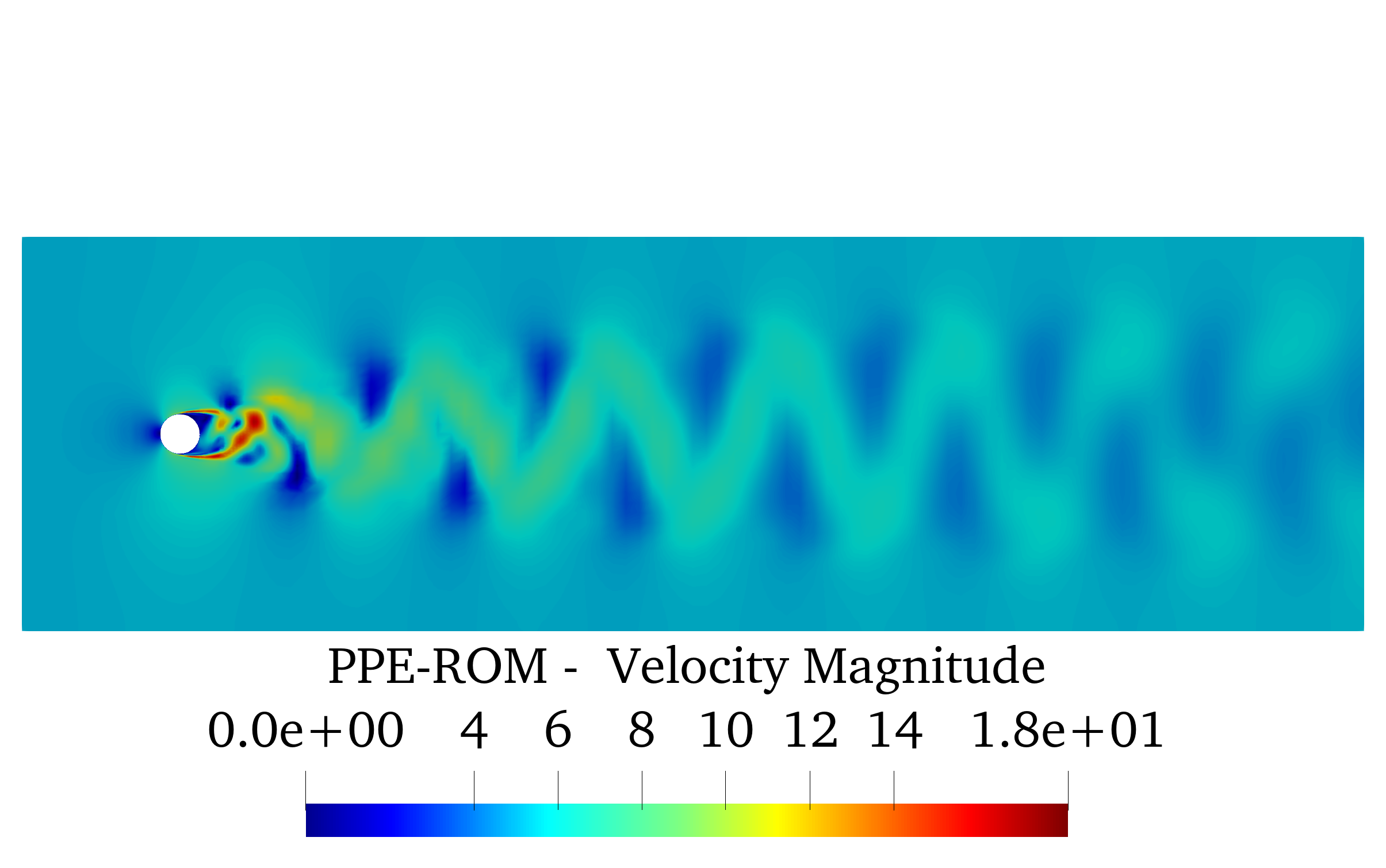}}
    \\ \vspace{-1.5cm}
    \subfloat[]{\includegraphics[width=0.46\textwidth]{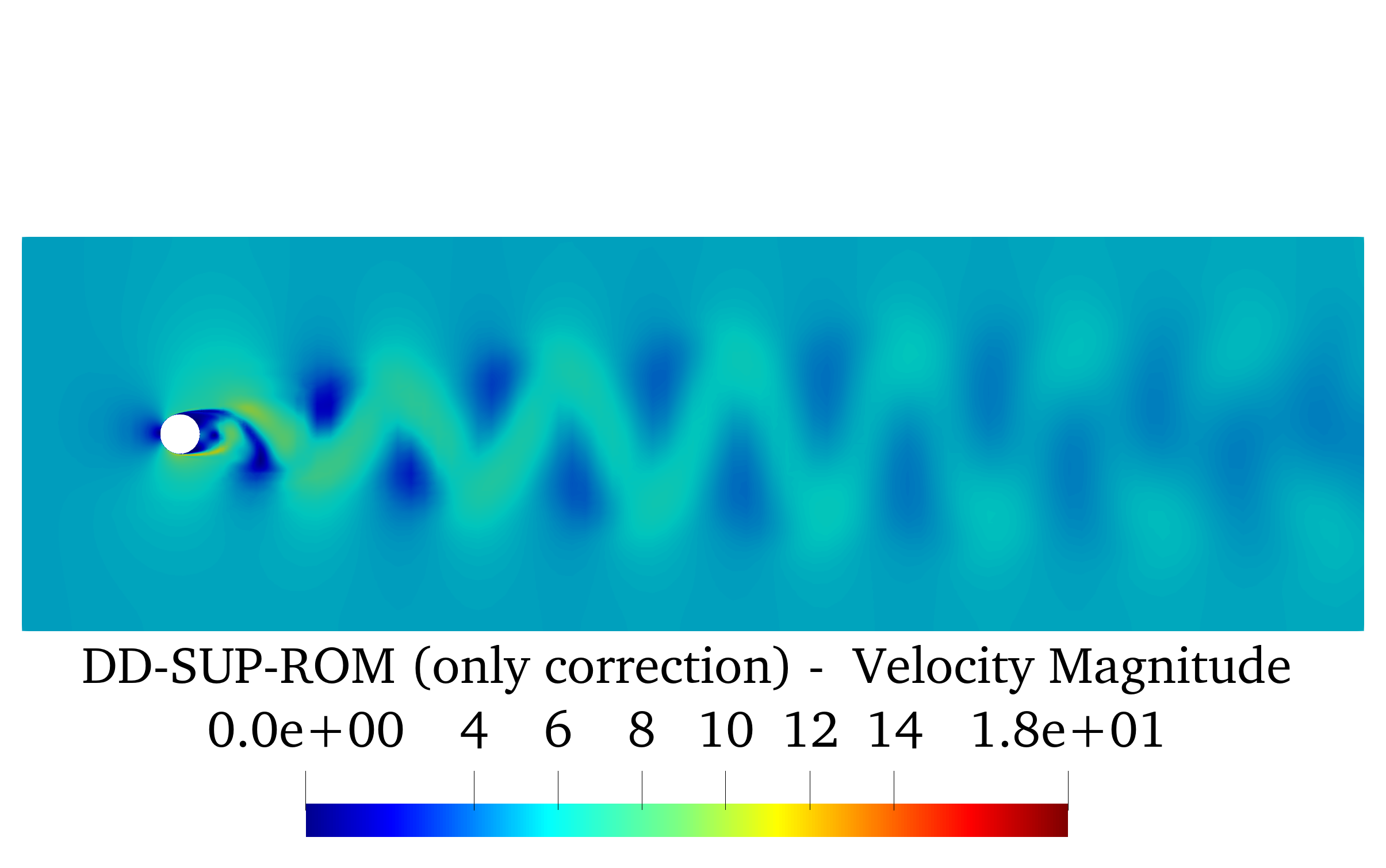}}
    \subfloat[]{\includegraphics[width=0.46\textwidth]{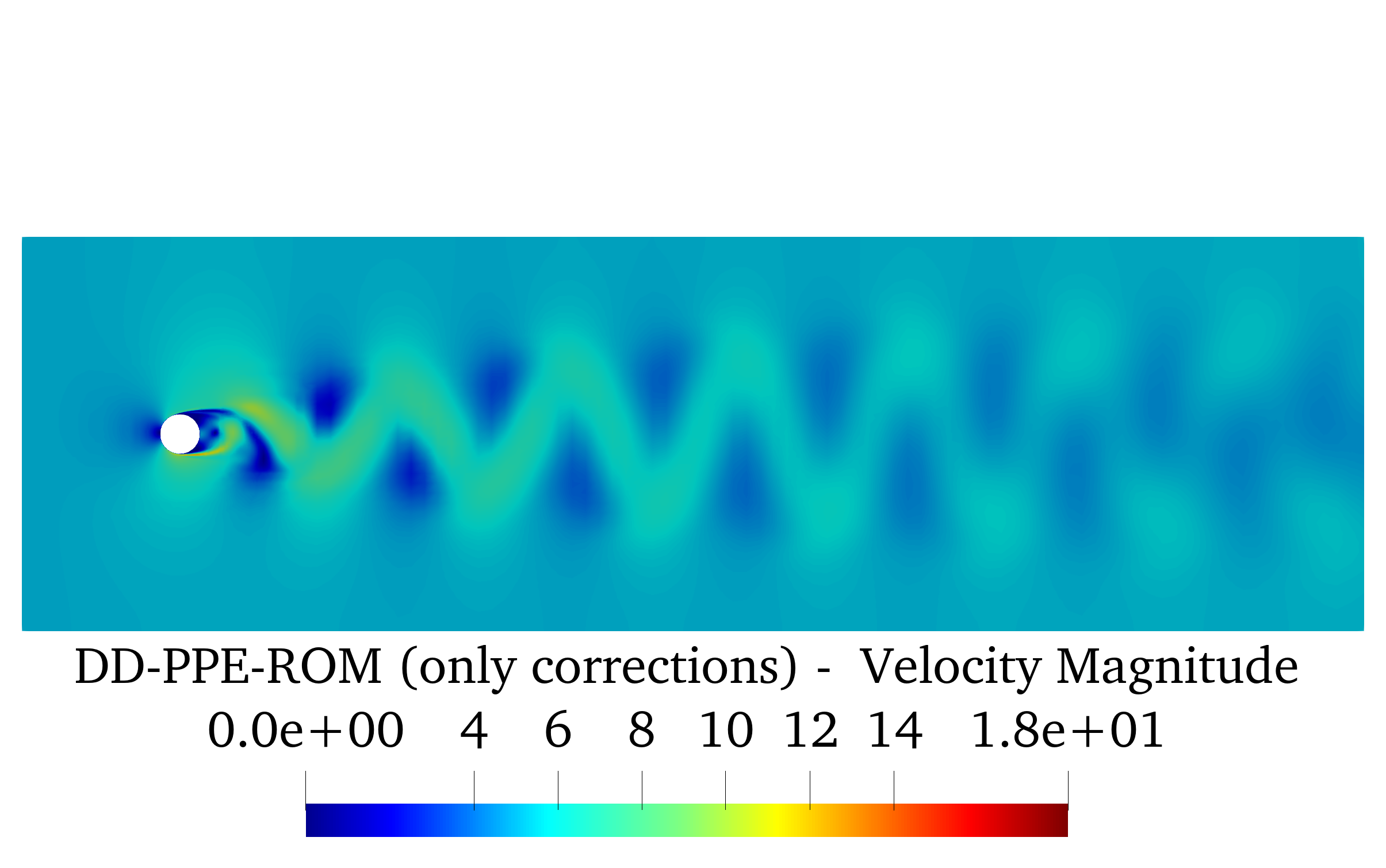}}\\ \vspace{-1.5cm}
    \subfloat[]{\includegraphics[width=0.46\textwidth]{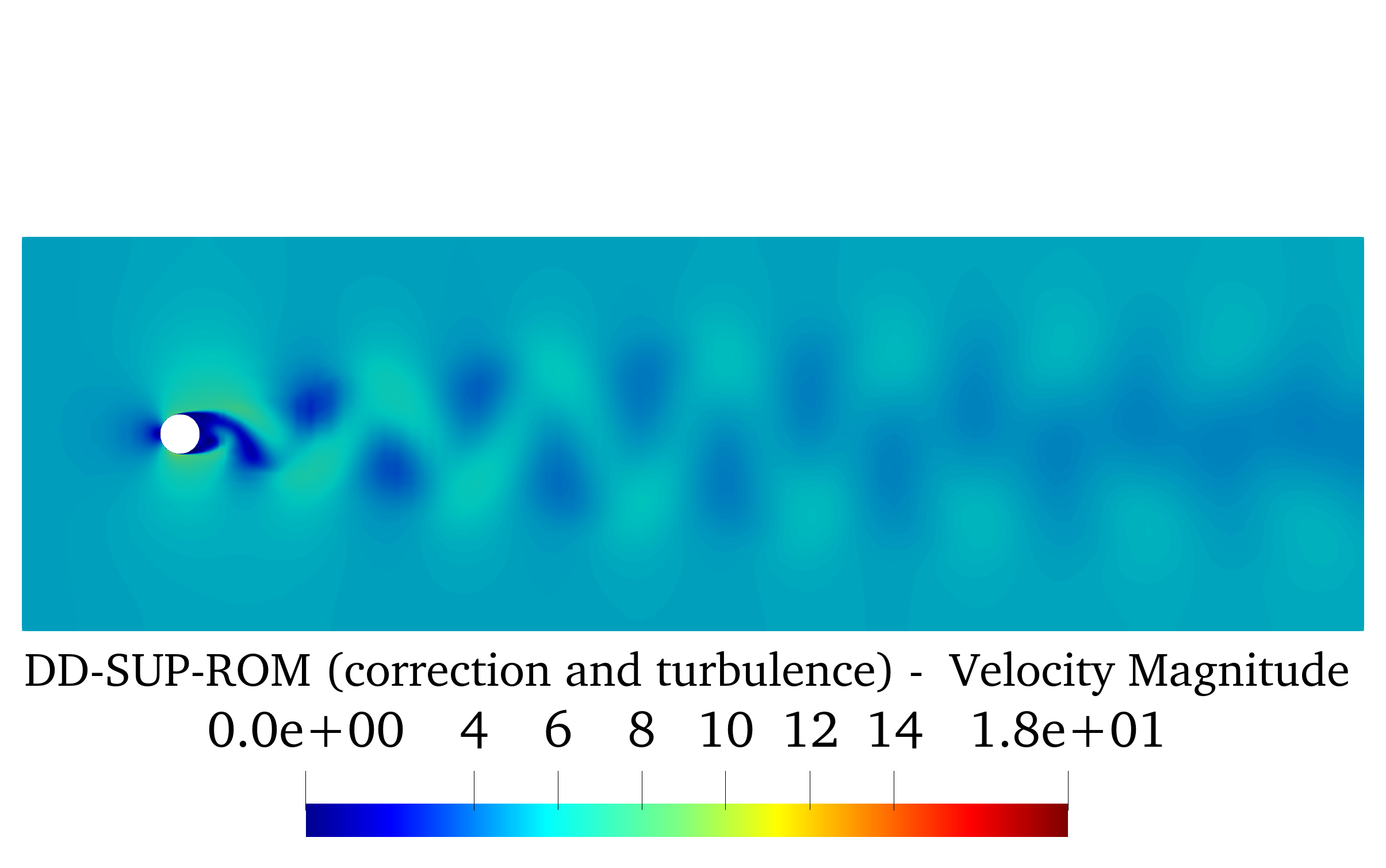}}
    \subfloat[]{\includegraphics[width=0.46\textwidth]{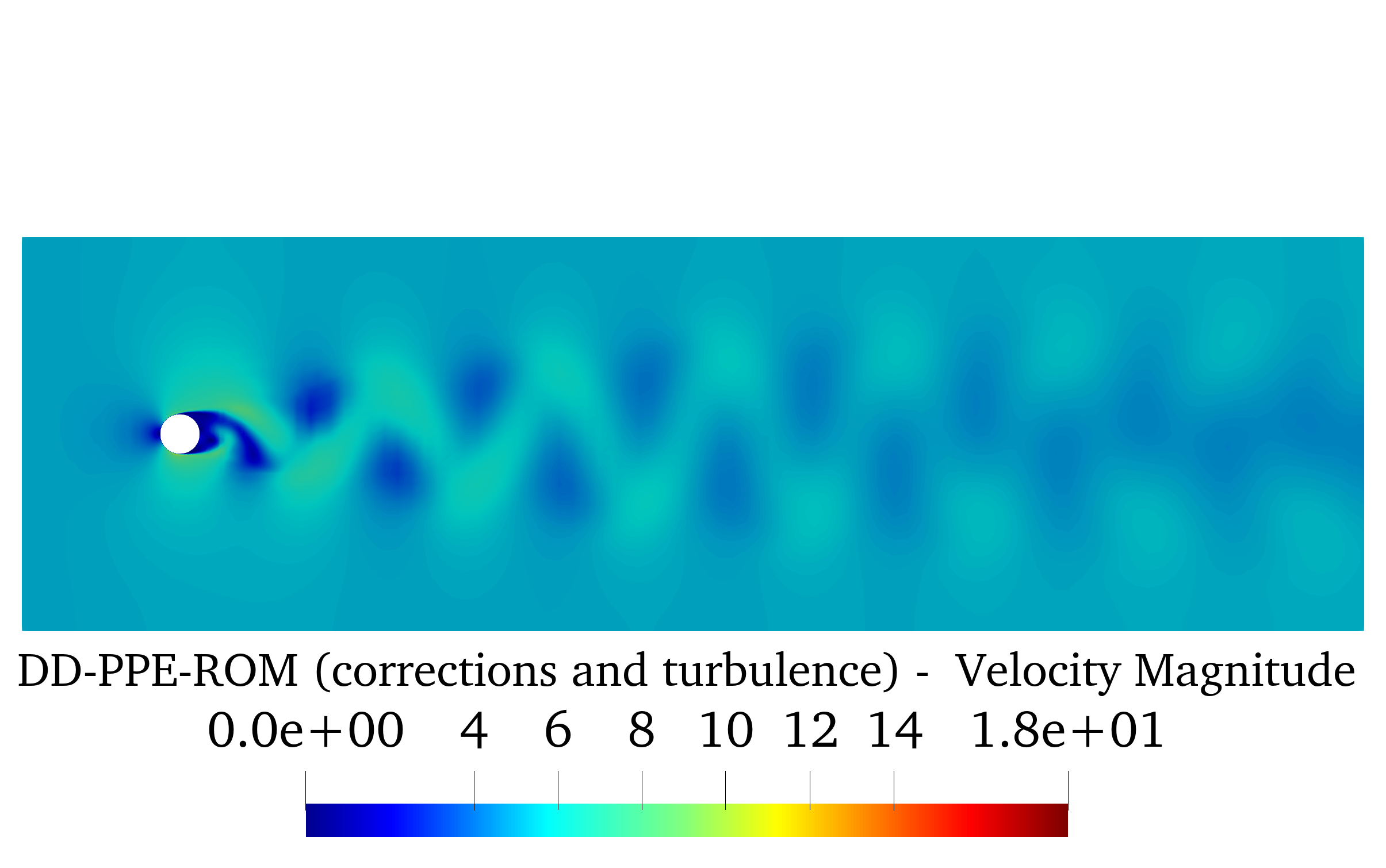}}\\ \vspace{-1.5cm}
    \subfloat[]{\includegraphics[width=0.46\textwidth]{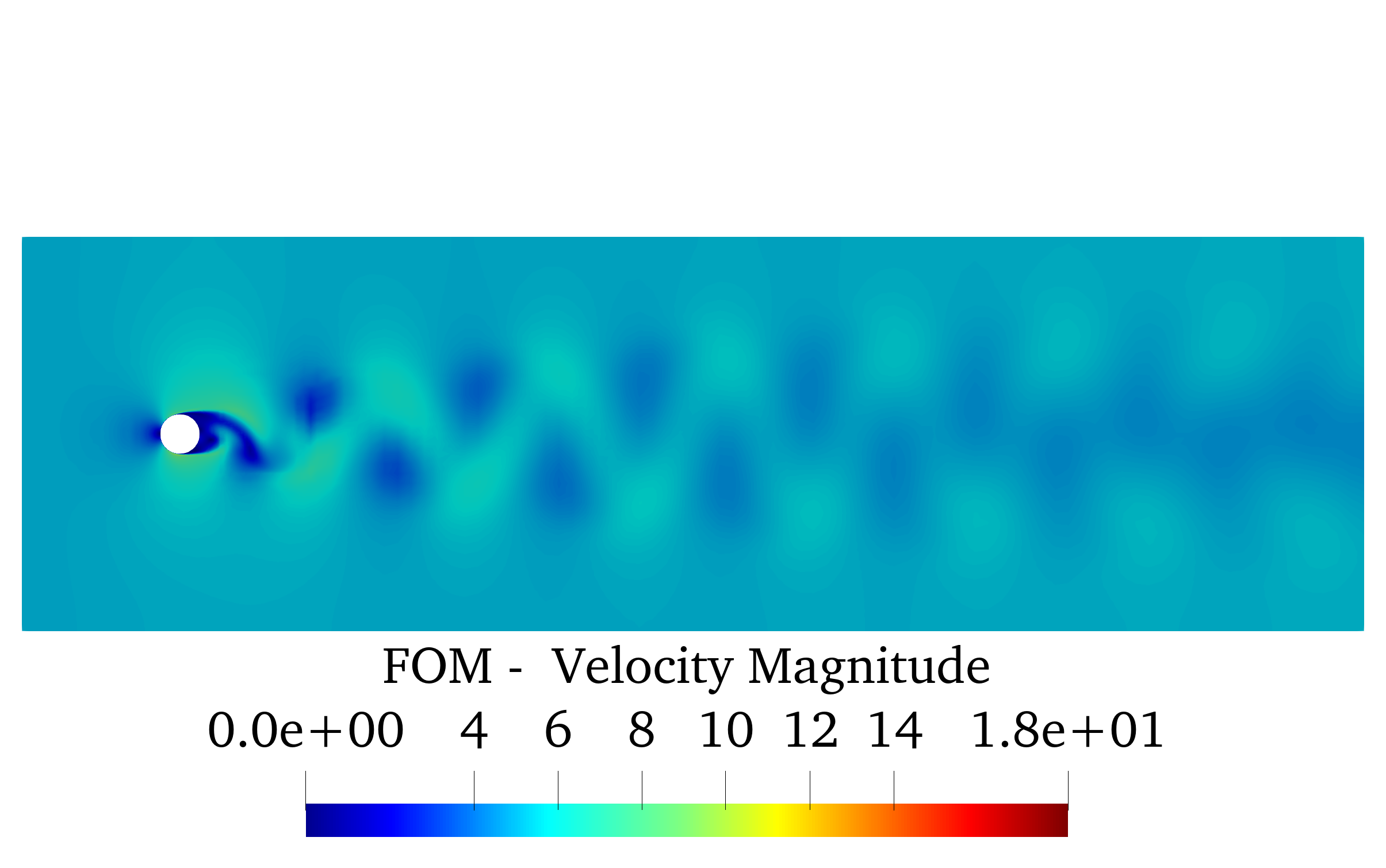}}
    \vspace{-0.4cm}
    \caption{Representation of the velocity magnitude field for the FOM, the SUP-ROM and the PPE-ROM simulations with and without the data-driven terms.}
    \label{paraview2}
\end{figure}

The pressure and the velocity magnitude fields are 
displayed in Figures \ref{paraview1} and \ref{paraview2}, respectively, for different SUP-ROM and PPE-ROM simulations.
The fields computed with the standard ROMs and those coming from the systems including the data-driven terms are different. 
In particular, the fields in panels (e) and (f) of Figures \ref{paraview1} and \ref{paraview2} 
are closer to the full order fields displayed in panel (g) of Figures \ref{paraview1} and \ref{paraview2}, especially in the region around the cylinder. The improvement of the 
accuracy nearby the circular cylinder is an important gain as it 
leads to a better reconstruction of the ROM lift coefficient.


\section{Conclusions and Outlook}
    \label{sec:conclusions}
One popular way to increase the accuracy of Galerkin ROMs in the under-resolved or marginally-resolved regimes is adding closure or correction terms
~\cite{ahmed2021closures}.
In this paper, the data-driven modeling of these correction terms, deeply analyzed in \cite{paper1}, 
was combined with 
turbulence modeling, using the eddy viscosity formulation developed in \cite{hijazi2020data}.
This yielded a novel hybrid data-driven closure strategy, which we tested numerically.
In our numerical investigation, we considered the two-dimensional flow past a circular cylinder at $Re= \num{50000}$ in the marginally-resolved regime.
We also considered several model configurations.

First, we tested two fundamentally different ROM pressure formulations:
(i) The supremizer ROM (SUP-ROM), in which additional (supremizer) modes for the velocity approximation are used in order to satisfy the inf-sup condition.
(ii) The pressure Poisson equation (PPE-ROM), in which the pressure approximation is determined by solving a Poisson equation instead of the continuity equation.

Secondly, a comparison between a first and a second order integration time scheme, in the resolution of the reduced system, 
was performed.
In our numerical investigation, we observed that the numerical dissipation associated with the second order scheme leads to a better reconstruction of the pressure and velocity fields when 
the novel data-driven terms are added. 
 
In conclusion, the 
hybrid data-driven closure strategy provides the best performance, leading to an increased accuracy of the reduced pressure and velocity field with respect to the high-fidelity solution, in both supremizer and Poisson approaches with a second order time integration scheme. 

\paragraph{Outlook}

There are several research directions that can be investigated:
\begin{itemize}
    \item The combined effect of correction terms and eddy viscosity modeling \green{is highly influenced by the nature of the test case taken into account. 
   The numerical method proposed} was tested in this paper only in the case of a 2D turbulent flow around a circular cylinder. Further tests with \green{unsteady problems, with }a more complex computational setting or with 
    3D grids \green{will be analysed in future works}.
    
    \green{Moreover, in this paper we focused on a single configuration for the training dataset: the snapshots are collected in a time interval of 20 seconds, the correction terms are built from data of 2 seconds and the online simulations last 8 seconds.
    In particular, the purely data-driven technique has a good extrapolation efficiency since an interval of 2 seconds (corresponding to 2 flow periods) collects enough information for the construction of the purely data-driven terms.
    
    Different settings would intuitively influence the performance of our numerical method and will be addressed in future works.}
    \item In this paper, the only parameter considered in the reduced order simulations is time. 
    Thus, the matrices appearing in the correction terms' ansatzes are parameter independent. An interesting task for the future would be the introduction of a parameter in the reduced formulation --- for instance the velocity at the inlet of the domain --- as in \cite{hijazi2020data}. In that case, the goal would be to express the data-driven terms as a function of both the parameter considered and 
    the time.
\end{itemize} 

\section{Acknowledgements}
We thank Prof. Claudio Canuto for his support.
We acknowledge the support by the European Commission H2020 ARIA (Accurate ROMs for Industrial Applications, GA 872442) project, by MIUR (Italian Ministry for Education University
and Research) and PRIN "Numerical Analysis for Full and Reduced
Order Methods for Partial Differential Equations" (NA-FROM-PDEs) project, and by the European
Research Council Consolidator Grant Advanced Reduced Order Methods with Applications in
Computational Fluid Dynamics-GA 681447, H2020-ERC COG 2015 AROMA-CFD. 
We also acknowledge support through National Science Foundation Grant Number DMS-2012253.
The main computations in this work were carried out by the usage of ITHACA-FV \cite{ithacasite}, a library maintained at
SISSA mathLab, an implementation in OpenFOAM \cite{ofsite} for reduced order modeling techniques. Its
developers and contributors are acknowledged.

\newpage

\bibliographystyle{plain}
\bibliography{main}

\end{document}